\documentclass[11pt,leqno]{article}
\normalfont
\renewcommand{\baselinestretch}{1.25}

\usepackage{graphicx}
\usepackage{amsfonts}
\usepackage{mathrsfs}
\usepackage{rotating}
\usepackage{mathdots}
\usepackage{mathabx}

\newcommand{\OverviewSection}{1}
  \newcommand{\IntroFibLucasTriangleFig}{Figure 1.1}
  \newcommand{\IntroGibTriangleFig}{Figure 1.2}

\newcommand{\RootSection}{2}
  \newcommand{\GMTWRootTheorem}{Theorem 2.1}
  \newcommand{\GibWProp}{Proposition 2.2}
  \newcommand{\GibRootTheorem}{Theorem 2.3}

\newcommand{\BinetSection}{3}
\newcommand{\ExplicitTheorem}{Theorem 3.1}
\newcommand{\ExplicitCorollaryOne}{Corollary 3.2}
\newcommand{\ExplicitCorollaryTwo}{Corollary 3.3}

\newcommand{\GameSection}{4}
  \newcommand{\GameFigureOne}{Figure 4.1}
  \newcommand{\GameFigureTwo}{Figure 4.2}
  \newcommand{\GameFigures}{Figures 4.1 and 4.2}
  \newcommand{\GameExample}{Example 4.1}
  \newcommand{\GameTheorem}{Theorem 4.2}
  \newcommand{\GameCorollary}{Corollary 4.3}
  \newcommand{\NumbersGame}{4}

\newcommand{\LatticeSection}{5}
\newcommand{ \RalphaGibFigure}{Figure 5.1}
  \newcommand{\InclusionExclusion}{Exercise 5.1}
  \newcommand{\RankedProp}{Proposition 5.2}
  \newcommand{\AlphaNProp}{Proposition 5.3}
  \newcommand{\RGFTheorem}{Theorem 5.4.A}
  \newcommand{\CountingTheorem}{Theorem 5.4.B}
  \newcommand{\EnumerativeTheorem}{Theorem 5.4}

\newcommand{\DiscussionSection}{6}

\newfont{\mysmallscbolditalics}{ecoc0500 at 10pt}
\newcommand{\mymathfrak}[1]{\mbox{\mysmallscbolditalics #1}}

\newfont{\myscbolditalics}{ecoc0500 at 11pt}

\newfont{\mybolditalics}{ecbi0500 at 11pt}

\newfont{\eulercursive}{eurm10 at 11pt}
\newcommand{\myq}{\mbox{\eulercursive q}}
\newcommand{\myr}{\mbox{\eulercursive r}}

\newcommand{\mya}{\mbox{\eulercursive a}}

\newcommand{\myl}{\mbox{\eulercursive `}}
\newfont{\smalleulercursive}{eurm10 at 9pt}

\newfont{\smallereulercursive}{eurm10 at 7pt}

\newcommand{\mysmallerq}{\mbox{\smallereulercursive q}}

\newfont{\myslantcyrillic}{wncyi10 at 11pt}

\newcommand{\QED}{\raisebox{0.5mm}{\fbox{\rule{0mm}{1.5mm}\ }}}

\newcounter{myfn}[page]
\renewcommand{\thefootnote}{\fnsymbol{footnote}}
\newcommand{\myfootnote}[1]{\setcounter{footnote}{\value{myfn}}%
    \footnote{#1}\stepcounter{myfn}}
    
\newcommand{\polyGoneone}[1]{{\overline{G}^{\mbox{\tiny \,$1,1$}}_{#1}}}
\newcommand{\polyGab}[1]{{\overline{G}^{\mbox{\tiny \,$\alpha,\beta$}}_{#1}}}

\newcommand{\polyGqone}[1]{{\overline{G}^{\mbox{\tiny \,$\mysmallerq,1$}}_{#1}}}
\newcommand{\polyGLucas}[1]{{\overline{G}^{\mbox{\tiny \,$2,1$}}_{#1}}}
\newcommand{\polyGFib}[1]{{\overline{G}^{\mbox{\tiny \,$1,1$}}_{#1}}}
\newcounter{rone}
\newcounter{rtwo}
\newcounter{rthree}
\newcounter{rfour}
\newcounter{rfive}
\newcounter{rsix}
\newcounter{rseven}
\newcommand{\RGF}{\mbox{\sffamily RGF}}

\newcommand{\CircleInteger}[1]{
\setlength{\unitlength}{0.14cm}
\begin{picture}(3,2) 
\put(1,1){\circle{2}}
\put(0.55,0.6){{\tiny #1}}
\end{picture}
}

\newcommand{\RightInterlace}{
\setlength{\unitlength}{0.5cm}
\begin{picture}(0.7,0) 
\thicklines
\multiput(0.05,0)(0.2,0){3}{\line(0,1){0.25}}
\multiput(0.15,0.15)(0.2,0){3}{\line(0,1){0.25}}
\end{picture}
}

\newcommand{\StrongInterlace}{
\setlength{\unitlength}{0.5cm}
\begin{picture}(0.7,0) 
\thicklines
\multiput(0.25,0)(0.2,0){2}{\line(0,1){0.25}}
\multiput(0.15,0.15)(0.2,0){3}{\line(0,1){0.25}}
\end{picture}
}

\newcommand{\Interlace}{
\setlength{\unitlength}{0.5cm}
\begin{picture}(0.9,0) 
\put(-0.4,0){
\setlength{\unitlength}{0.5cm}
\begin{picture}(0,0) 
\put(-0.125,0){\tiny (}
\put(0.1,0){\tiny )}
\thicklines
\put(0.095,0){\line(0,1){0.25}}
\multiput(0.4,0)(0.2,0){2}{\line(0,1){0.25}}
\multiput(0.3,0.15)(0.2,0){3}{\line(0,1){0.25}}
\end{picture}}
\end{picture}
}

\newcommand{\VertexBlankFib}{
\setlength{\unitlength}{1.5cm}
\begin{picture}(0,0)
\put(-0.25,0){
\begin{picture}(0,0)
\put(0,0){\circle*{0.1}} 
\end{picture}
}
\end{picture}
}

\begin{document}
\pagenumbering{arabic}
\thispagestyle{empty}%
\vspace*{-0.7in}
\hfill \parbox{2.54in}{\hfill {\scriptsize December 29, 2020}} 

\begin{center}
{\large \bf Sign-alternating Gibonacci Polynomials} 

\vspace*{0.05in}
\renewcommand{\thefootnote}{1}
Robert G.\ Donnelly\footnote{Department of Mathematics and Statistics, Murray State
University, Murray, KY 42071\\ 
\hspace*{0.25in}Email: {\tt rob.donnelly@murraystate.edu}}, 
\renewcommand{\thefootnote}{2} 
\hspace*{-0.07in}Molly W.\ Dunkum\footnote{Department of Mathematics, Western Kentucky University, Bowling Green, KY 42101\\ 
\hspace*{0.25in}Email: {\tt molly.dunkum@wku.edu}}, 
\renewcommand{\thefootnote}{3} 
\hspace*{-0.07in}Murray L. Huber\footnote{Department of Mathematics, Western Kentucky University, Bowling Green, KY 42101\\ 
\hspace*{0.25in}Email: {\tt murraylhuber@gmail.com}}, 
\renewcommand{\thefootnote}{4} 
\hspace*{-0.07in}and Lee Knupp\footnote{Department of Mathematics, Western Kentucky University, Bowling Green, KY 42101\\ 
\hspace*{0.25in}Email: {\tt knupp.19@osu.edu}}  

\end{center} 

\begin{abstract}
We consider various properties and manifestations of some sign-alternating univariate polynomials borne of right-triangular integer arrays related to certain generalizations of the Fibonacci sequence. 
Using a theory of the root geometry of polynomial sequences developed by J.\ L.\ Gross, T.\ Mansour, T.\ W.\ Tucker, and D.\ G.\ L.\ Wang, we show that the roots of these `sign-alternating Gibonacci polynomials' are real and distinct, and we obtain explicit bounds on these roots. 
We also derive Binet-type closed expressions for the polynomials. 
Some of these results are applied to resolve finiteness questions pertaining to a one-player combinatorial game (or puzzle) modelled after a well-known puzzle we call the `Networked-numbers Game.' 

Elsewhere, the first- and second-named authors, in collaboration with A.\ Nance, have found rank symmetric `diamond-colored' distributive lattices naturally related to certain representations of the special linear Lie algebras. 
Those lattice cardinalities can be computed using sign-alternating Fibonacci polynomials, and the lattice rank generating functions correspond to the rows of some new and easily defined triangular integer arrays. 
Here, we present Gibonaccian, and in particular Lucasian, versions of those symmetric Fibonaccian lattices/results, but without the algebraic context of the latter. 

\ 

\begin{center}
\renewcommand{\thefootnote}{$\mbox{\ }$}
\footnote{\ \ }
\end{center} 
\end{abstract}

\noindent 
{\large \bf \S \OverviewSection.\ Introduction.} 
The sign-alternating polynomials we consider here naturally generalize a variation of the so-called Fibonacci polynomials defined in \S 37 of \cite{Koshy} and \S 9.4 of \cite{BQ}.  
Coefficients for the latter polynomials can be viewed as a nicely-structured right-triangular array of positive integers --- see the leftmost triangle of \IntroFibLucasTriangleFig\ below or OEIS A011973 \cite{OEIS}. 
To introduce the sign-alternating Gibonacci\myfootnote{The adjective ``Gibonacci'' is a portmanteau identifying certain {\bf G}eneralized f $\!\!${\bf IBONACCI} sequences. 
In \S 7 of \cite{Koshy}, T.\ Koshy attributes this neologism to A.\ T.\ Benjamin and J.\ J.\ Quinn from their well-known book \cite{BQ}.} polynomials, we consider a two-parameter generalization of the foregoing right-triangular Fibonacci array, cf.\ \IntroGibTriangleFig.  
Throughout this paper, these parameters are positive real numbers $\alpha$ and $\beta$, although, for enumerative reasons, at times we take special interest in those cases where $\alpha$ and $\beta$ are integers.
We note here at the outset that variations of our sign-alternating Gibonacci polynomials occur as special cases of the so-called Horadam sequence of polynomials defined in \cite{HK} and of the type $(0,1)$ polynomials of \cite{GMTW1} . 

\begin{figure}[hb]
\begin{center}
{\renewcommand{\baselinestretch}{1.2}
\footnotesize
{\begin{tabular}{ccccc}
1 & & & & \\
1 & & & & \\
1 & 1 & & & \\
1 & 2 & & & \\
1 & 3 & 1 & & \\
1 & 4 & 3 & & \\
1 & 5 & 6 & 1 & \\
1 & 6 & 10 & 4 & \\
1 & 7 & 15 & 10 & 1 \\
1 & 8 & 21 & 20 & 5 \\
 & & $\cdots$ & & 
\end{tabular}}
\hspace*{1.5in}
{\begin{tabular}{ccccc}
2 & & & & \\
1 & & & & \\
1 & 2 & & & \\
1 & 3 & & & \\
1 & 4 & 2 & & \\
1 & 5 & 5 & & \\
1 & 6 & 9 & 2 & \\
1 & 7 & 14 & 7 & \\
1 & 8 & 20 & 16 & 2 \\
1 & 9 & 27 & 30 & 9 \\
& & $\cdots$ & & 
\end{tabular}}
}

\vspace*{0.1in}
{\small \IntroFibLucasTriangleFig: The right-triangular Fibonacci (left) and Lucas (right) arrays}
\end{center}
\end{figure}

What we call the {\em $(\alpha,\beta)$-Gibonacci right-triangular array} (or {\em Gibonacci array}, for short) is the array 
\[\mymathfrak{G}(\alpha,\beta) := (g^{\mbox{\tiny $\alpha,\beta$}}_{k,j})_{k\in\{0,1,2,\ldots\}, j\in\{0,1,\ldots,\lfloor\frac{k}{2}\rfloor\}}\] 
where $g^{\mbox{\tiny $\alpha,\beta$}}_{0,0} := \alpha$, $g^{\mbox{\tiny $\alpha,\beta$}}_{1,0} := \beta$, and 
\[g^{\mbox{\tiny $\alpha,\beta$}}_{k,j} := g^{\mbox{\tiny $\alpha,\beta$}}_{k-1,j} + g^{\mbox{\tiny $\alpha,\beta$}}_{k-2,j-1}\] 
for integers $j$ and $k$ with $k \geq 2$ and with the understanding that $g^{\mbox{\tiny $\alpha,\beta$}}_{k,j} := 0$ when $j < 0$ or $j > \lfloor \frac{k}{2} \rfloor$. 
The first several rows of the array are depicted in \IntroGibTriangleFig.   

\begin{figure}[ht]
\begin{center}
{\renewcommand{\baselinestretch}{1.2}
\footnotesize
\begin{tabular}{ccccc}
$\alpha$ & & & & \\
$\beta$ & & & & \\
$\beta$ & $\alpha$ & & & \\
$\beta$ & $\alpha+\beta$ & & & \\
$\beta$ & $\alpha+2\beta$ & $\alpha$ & & \\
$\beta$ & $\alpha+3\beta$ & $2\alpha+\beta$ & & \\
$\beta$ & $\alpha+4\beta$ & $3\alpha+3\beta$ & $\alpha$ & \\
$\beta$ & $\alpha+5\beta$ & $4\alpha+6\beta$ & $3\alpha+\beta$ & \\
$\beta$ & $\alpha+6\beta$ & $5\alpha+10\beta$ & $6\alpha+4\beta$ & $\alpha$ \\
$\beta$ & $\alpha+7\beta$ & $6\alpha+15\beta$ & $10\alpha+10\beta$ & $4\alpha+\beta$ \\
& & $\cdots$ & & 
\end{tabular}}

\vspace*{0.1in}
{\small \IntroGibTriangleFig: The right-triangular Gibonacci array $\mymathfrak{G}(\alpha,\beta)$}
\end{center}
\end{figure}

The right-triangular arrays $\mymathfrak{G}(1,1)$ and $\mymathfrak{G}(2,1)$ are partially depicted in \IntroFibLucasTriangleFig. 
Row sums of $\mymathfrak{G}(1,1)$ correspond to the Fibonacci sequence $f_{0}=1$, $f_{1}=1$, $f_{2}=2$, $f_{3}=3$, $f_{4}=5$, etc, while row sums of $\mymathfrak{G}(2,1)$ correspond to the Lucas sequence $\myl_{0} = 2$, $\myl_{1} = 1$, $\myl_{2}=3$, $\myl_{3}=4$, $\myl_{4}=7$, etc. 
So $\mymathfrak{G}(1,1)$ (respectively, $\mymathfrak{G}(2,1)$) is the right-triangular Fibonacci (resp.\ Lucas) array. 

For the remainder of this section, let $k$ and $j$ be integers with $0 \leq j \leq \lfloor k/2 \rfloor$. 
Using induction and the defining recurrence for $\mymathfrak{G}(\alpha,\beta)$, it is easy to see that 
\[g^{\mbox{\tiny $\alpha,\beta$}}_{k,j} = {k-j-1 \choose j-1}\alpha + {k-j-1 \choose j}\beta,\]
with the usual convention that for integers $a$ and $b$, the binomial coefficient ${a \choose b}$ is $0$ when we do not have $0 \leq b \leq a$. 
When $\alpha$ and $\beta$ are positive integers, the number $g^{\mbox{\tiny $\alpha,\beta$}}_{k,j}$ has a nice combinatorial interpretation as the number of ``$(\alpha,\beta)$-phased $k$-tilings,'' cf.\ Combinatorial Theorem 13 from \cite{BQ}. 

The focus of this paper is on the following polynomials: \[\polyGab{k}(x) := \sum_{j=0}^{\lfloor \frac{k}{2} \rfloor} (-1)^{j}g^{\mbox{\tiny $\alpha,\beta$}}_{k,j}\, x^{\lfloor \frac{k}{2} \rfloor - j},\] 
so $\polyGab{n}(x)$ is a polynomial whose integer coefficients are signed versions of the $k^{\mbox{\tiny th}}$ row of $\mymathfrak{G}(\alpha,\beta)$.  
(The ``{\verb \overline{G} }'' LaTeX notation used above is meant to be a visual reminder that the signs of the coefficients alternate.)
As an example, when $\alpha = 2$ and $\beta = 1$, then $\overline{G}^{\mbox{\tiny 2,1}}_{7}(x) = x^{3} - 7x^{2} + 14x - 7$.  
We call these the {\em sign-alternating $(\alpha,\beta)$-Gibonacci array polynomials}, or for brevity {\em sign-alternating Gibonacci polynomials}.  
When needed, we set $\polyGab{-1}(x) := 0$. 
From the defining recurrence of the Gibonacci array we get the following fundamental recurrence for our sign-alternating Gibonacci polynomials:  
\[ \polyGab{k}(x) = x^{(k-1)\, \mbox{\scriptsize mod 2}} \polyGab{k-1}(x) - \polyGab{k-2}(x)\]
for all integers $k \geq 2$, where $\polyGab{0}(x) = \alpha$ and $\polyGab{1}(x) = \beta$.  
A routine application of the preceding recurrence is the following, for all $k \geq 2$: 
\[ \polyGab{k}(x) = x^{(k-1)\, \mbox{\scriptsize mod 2}} \beta\, \polyGoneone{k-1}(x) - \alpha\, \polyGoneone{k-2}(x),\]
where any $\polyGoneone{m}(x)$ is to be viewed as a sign-alternating Fibonacci polynomial obtained from the Fibonacci array $\mymathfrak{G}(1,1)$. 

We have encountered certain of these polynomials on two separate occasions within the context of our work in combinatorial representation theory: once in our study of a game of numbers related to Weyl group actions on Weyl symmetric functions, and once in our study of distributive lattice models for semisimple Lie algebra representations. 
A primary purpose of this paper is to draw attention to the titular family of sign-alternating Gibonacci polynomials by exhibiting aspects of these disparate connections.   

The first of the two aforementioned connections is to a single-player game (or puzzle) most often called the `Numbers Game', see \cite{E} or \S 4.3 of \cite{BB}. 
We prefer to call this puzzle the `Networked-numbers Game' to emphasize the crucial role of an underlying graph that links the numbers together.   
In analyzing a two-node version of this game, we rediscovered certain numerical constraints on game play that were originally found by Eriksson \cite{E}. 
Our new proof of this result is obtained by directly relating said constraints to the roots of the sign-alternating Fibonacci polynomials. 
This result is recovered in \S \GameSection\ here as a corollary (\GameCorollary) of a more general result (\GameTheorem). 
The latter result concerns roots of sign-alternating Gibonacci polynomials, and motivates the considerations of \S \RootSection\ and \S \BinetSection\ here. 

The second connection is to some distributive lattices that have many manifestations in the literature but were re-discovered by us in a Lie representation theoretic context.  
Specifically, in \cite{DDMN} are presented some symmetric distributive lattices related to the sequence of `symmetric' Fibonacci numbers $1$, $3$, $8$, $21$, $55$, etc. 
The `symmetric Fibonaccian lattices' of that paper are shown to be models for certain representations of the special linear Lie algebras and for the related skew Schur functions. 
Some new enumerative identities relating to those lattices were also obtained in \cite{DDMN}.
Those results motivated us to find Gibonaccian ranked poset analogs of the symmetric Fibonaccian lattices, which we present in \S \LatticeSection\ here. 
These new families of ranked posets include what we call symmetric Lucasian lattices, which are symmetric distributive lattices related to the sequence of `symmetric' Lucas numbers $2$, $3$, $7$, $18$, $47$, etc.  
In general, our new posets possess enumerative properties (\EnumerativeTheorem) that analogize results for the symmetric Fibonaccian lattices but without (as far as we can tell) an analogous algebraic context.  
Our work in this section generalizes some of the results of \cite{Munarini--Salvi} which, in our notation, was mostly concerned with the sequence of rank sizes for the $n=3$ and $(\alpha,\beta) \in \{(1,1),(2,1)\}$ cases.

\vspace*{0.1in}
\noindent 
{\large \bf \S \RootSection.\ Roots of sign-alternating Gibonacci polynomials.} 
A systematic study of the root geometry of certain recursively defined polynomial sequences is undertaken in \cite{GMTW1} and \cite{GMTW2}. 
Here, we connect our sign-alternating Gibonacci polynomials to the environment of  \cite{GMTW1} so we can use the framework of that paper to obtain many of the root-related results of this section. 
Some terminology: For sets of real numbers $X = \{x_{1},\ldots,x_{n}\}$ and $Y = \{y_{1},\ldots,y_{m}\}$ each indexed from smallest to largest, say $X$ {\em interlaces} $Y$ {\em from both sides} and write $X \StrongInterlace Y$ if $m=n-1$ and 
$x_{1} < y_{1} < x_{2} < y_{2} < \cdots < y_{n-1} < x_{n}$;  
say $X$ {\em interlaces} $Y$ {\em from the right} and write $X \RightInterlace Y$ if $m=n$ and
$y_{1} < x_{1} < y_{2} < x_{2} < \cdots < y_{n} < x_{n}$; and say $X$ {\em interlaces} $Y$ and write  $X \Interlace Y$ if either $X \StrongInterlace Y$ or $X \RightInterlace Y$. 
If a sequence $\{x_{k}\}$ converges to some $x$ strictly-increasingly (resp.\ strictly-decreasingly), we write $\{x_{k}\} \nearrow x$ (resp.\ $\{x_{k}\} \searrow x$).
Throughout this section, we set $\myq := \alpha/\beta$. 
Since $\mbox{sgn}\left(\rule[-2.5mm]{0mm}{5.5mm}\polyGab{k}(x)\right) = \mbox{sgn}\left(\rule[-2.5mm]{0mm}{5.5mm}\polyGqone{k}(x)\right)$ for any real number $x$, then  some results of this section will be developed for $\polyGqone{k}(x)$ rather than $\polyGab{k}(x)$. 

We now develop a special case of Theorem 2.6 of \cite{GMTW1}, from which we will deduce most of our root-related results for sign-alternating Gibonacci polynomials $\polyGqone{k}(x)$. 
Define a sequence $\{W_{k}(x)\}_{k \geq 0}$ of polynomials recursively by the rule $W_{k}(x) = W_{k-1}(x)+xW_{k-2}(x)$ when $k \geq 2$, with $W_{0}(x) := 1$ and $W_{1}(x) := 1+\myq x$. 
Let $d'_{k}$ be the number of real roots of $W_{k}(x)$, and let $R_{k} := \{\xi_{k,1},\ldots,\xi_{k,d'_{k}}\}$ be the set of these roots, indexed from smallest to largest. 

\noindent 
{\bf \GMTWRootTheorem}\ \ {\sl With $\{W_{k}(x)\}_{k \geq 0}$ defined as above, then $W_{k}(x)$ is a polynomial of degree $d'_{k} = \lfloor (k+1)/2 \rfloor$ and has $d'_{k}$ distinct real roots. 
For any $i \geq 1$, the sequence $\{\xi_{k,i}\}_{k \geq 1}$ converges to $-\infty$. 
Also, for any $k \geq 1$}, $R_{k+1} \Interlace R_{k}$ {\sl and $R_{k+2} \StrongInterlace R_{k}$.}  
{\sl If $\myq \leq 2$, then for any $i \geq 0$, $\{\xi_{k,d'_{k}-i}\}_{k \geq 1} \nearrow -\frac{1}{4}$.} 
{\sl If $\myq > 2$, then $\{\xi_{k,d'_{k}}\}_{k \geq 1} \nearrow \frac{-\mysmallerq+1}{\mysmallerq^{2}}$ and, for any $i \geq 1$, $\{\xi_{k,d'_{k}-i}\}_{k \geq 1} \nearrow -\frac{1}{4}$.} 

{\em Proof.} This is a special case of Theorem 2.6 of [GMTW1], where we have taken $a=b=1$, $c=0$, $r = -1/\myq$, $t = \myq$, $x^{*} = -1/4$, $r^{*} = -\frac{1}{4}-\frac{1}{2\mysmallerq}$, and $\displaystyle y^{*} = \left\{\begin{array}{cc}\frac{-\mysmallerq+1}{\mysmallerq^{2}} & \mbox{ if } \myq \geq 1\\ 0 & \mbox{ if } \myq < 1\end{array}\right.$.\hfill\QED  

\noindent
{\bf \GibWProp}\ \ {\sl We have} $\polyGqone{k}(x) = x^{\lfloor k/2 \rfloor}W_{k-1}(-\frac{1}{x})$ {\sl for any $x \ne 0$ and for any $k \geq 2$.}

{\em Proof.} Observe that the expression $x^{\lfloor k/2 \rfloor}W_{k-1}(-\frac{1}{x})$ is defined for all $x \ne 0$, has a removable discontinuity at $x=0$, and simplifies to a polynomial. 
So we let $\overline{P}_{k}(x)$ be the polynomial simplification of $x^{\lfloor k/2 \rfloor}W_{k-1}(-\frac{1}{x})$. 
It is routine to verify that the $\overline{P}_{k}(x)$'s satisfy the same recurrence relations as the $\polyGqone{k}(x)$'s with the same initial conditions.\hfill\QED

\noindent
{\bf \GibRootTheorem}\ \ {\sl For any $k \geq 2$, the degree $d_{k} := \lfloor k/2 \rfloor$ polynomial $\polyGab{k}(x)$ has $d_{k}$ distinct positive real roots which we gather in the set $S_{k} := \{\zeta_{k,1},\ldots,\zeta_{k,d_{k}}\}$, indexed from smallest to largest. 
For any $i \geq 1$, the sequence $\{\zeta_{k,i}\}_{k \geq 2}$ converges to $0$. 
Also, for any $k \geq 2$}, $S_{k+1} \Interlace S_{k}$ {\sl and $S_{k+2} \StrongInterlace S_{k}$.}  
{\sl If $\alpha/\beta \leq 2$, then for any $i \geq 0$, $\{\zeta_{k,d_{k}-i}\}_{k \geq 2} \nearrow 4$.} 
{\sl Now suppose $\alpha/\beta > 2$.  
Then $\{\zeta_{k,d_{k}}\}_{k \geq 2} \nearrow \frac{(\alpha/\beta)^{2}}{\alpha/\beta - 1}$ and, for any $i \geq 1$, $\{\zeta_{k,d_{k}-i}\}_{k \geq 2} \nearrow 4$.} 

{\em Proof.} All claims follow by putting \GibWProp\ together with \GMTWRootTheorem.\hfill\QED

\vspace*{0.1in}
\noindent 
{\large \bf \S \BinetSection. Binet-type formulas for sign-alternating Gibonacci polynomials.} 
Here we show how to use a rudimentary linear algebraic approach 
to derive Binet-type formulas for the sign-alternating Gibonacci polynomials. 
Such methodology is standard in solving recurrences such as those that define our sign-alternating Gibonacci polynomials (see, for example, \cite{Benoumhani}, \cite{HK}). 
For any positive integer $m$, the fundamental recurrence for these polynomials yields the matrix identities 
\begin{eqnarray*}
\left(\begin{array}{c} \polyGab{2m}(x)\\ \polyGab{2m+1}(x)\end{array}\right) & \mbox{\huge $=$} & 
\left(\begin{array}{cc} -1 & x\\ -1 & x-1\end{array}\right)\left(\begin{array}{c} \polyGab{2m-2}(x)\\ \polyGab{2m-1}(x)\end{array}\right)\\
& \mbox{\huge $=$} & \rule[-3mm]{0mm}{9mm}\mbox{\huge $\cdots$}\\ 
& \mbox{\huge $=$} & \left(\begin{array}{cc} -1 & x\\ -1 & x-1\end{array}\right)^{m}\left(\begin{array}{c} \polyGab{0}(x)\\ \polyGab{1}(x)\end{array}\right)\\ 
& \mbox{\huge $=$} & \left(\begin{array}{cc} -1 & x\\ -1 & x-1\end{array}\right)^{m}\left(\begin{array}{c} \alpha\\ \beta\end{array}\right).
\end{eqnarray*}
One can check that the transition matrix $\left(\begin{array}{cc} -1 & x\\ -1 & x-1\end{array}\right)$ has eigenvalues 
\[\lambda = \frac{1}{2}\left(x-2+\sqrt{x^{2}-4x}\right) \hspace*{0.25in}\mbox{and}\hspace*{0.25in} \kappa = \frac{1}{2}\left(x-2-\sqrt{x^{2}-4x}\right)\] with corresponding eigenvectors $\left(\begin{array}{c}x-1-\lambda\\ 1\end{array}\right)$ and $\left(\begin{array}{c}x-1-\kappa\\ 1\end{array}\right)$ respectively. 
The eigenvalues $\lambda$ and $\kappa$ are roots of the characteristic polynomial \[\det\left(\begin{array}{cc} t+1 & -x\\ 1 & t-(x-1)\end{array}\right) = t^{2}-(x-2)t+1\] of our transition matrix, so $\lambda\kappa = 1$ and $\lambda+\kappa = x-2$. 
Using these quantities to diagonalize the transition matrix yields the following expression for its $m^{\mbox{\tiny th}}$ power: 
\begin{center}
{\small 
$\displaystyle \left(\begin{array}{cc} -1 & x\\ -1 & x-1\end{array}\right)^{m} =\rule[0mm]{130mm}{0mm}$\\ 
$\displaystyle \frac{1}{\kappa-\lambda}\left(\begin{array}{cc} \lambda^{m}(x-1-\lambda)-\kappa^{m}(x-1-\kappa) & -\lambda^{m}(x-1-\lambda)(x-1-\kappa)+\kappa^{m}(x-1-\kappa)(x-1-\lambda)\\ \lambda^{m}-\kappa^{m} & -\lambda^{m}(x-1-\kappa)+\kappa^{m}(x-1-\lambda)\end{array}\right)$. 
} 
\end{center}
Then 
$\displaystyle \left(\begin{array}{c} \polyGab{2m}(x)\\ \polyGab{2m+1}(x)\end{array}\right) = \left(\begin{array}{cc} -1 & x\\ -1 & x-1\end{array}\right)^{m}\left(\begin{array}{c} \alpha\\ \beta\end{array}\right) = $
\begin{center}
{\small 
$\displaystyle \frac{1}{\kappa-\lambda}\left(\begin{array}{c} \lambda^{m}(x-1-\lambda)[\alpha-(x-1-\kappa)\beta]-\kappa^{m}(x-1-\kappa)[\alpha-(x-1-\lambda)\beta] \\ \lambda^{m}[\alpha-(x-1-\kappa)\beta]-\kappa^{m}[\alpha-(x-1-\lambda)\beta]\end{array}\right)$. 
} 
\end{center}
We therefore obtain nonrecursive closed-form expressions for $\polyGab{2m}(x)$ and $\polyGab{2m+1}(x)$. 
We formally record these results in the following slightly modified forms. 

\noindent 
{\bf \ExplicitTheorem}\ \ {\sl Keep the above notation.  For any nonnegative integer $m$, we have}
\setlength{\unitlength}{1in}
\begin{picture}(0,0)
\put(2.55,-0.9){\QED}
\end{picture}
\begin{eqnarray*}
\polyGab{2m}(x) & = & \frac{\lambda^{m}[(\kappa+1)\alpha-x\beta]-\kappa^{m}[(\lambda+1)\alpha - x\beta]}{\kappa-\lambda}\\ 
\polyGab{2m+1}(x) & = & \frac{\lambda^{m}[\alpha-(\lambda+1)\beta]-\kappa^{m}[\alpha - (\kappa+1)\beta]}{\kappa-\lambda}
\end{eqnarray*}

As an illustration of this theorem, we note that the preceding formulas nicely specialize for the $\alpha= 1,\beta = 1$ (Fibonacci) and $\alpha=2,\beta=1$ (Lucas) cases. 
What follows next are just variations on well-known closed-form expressions for related Fibonacci/Lucas polynomials. 

\noindent
{\bf \ExplicitCorollaryOne}\ \ {\sl Let $m$ be a nonnegative integer.  The sign-alternating Fibonacci polynomials can be written thusly:} 
\[\polyGFib{2m}(x) = \frac{\lambda^{2m+1}-1}{\lambda^{m}(\lambda-1)} \hspace*{0.25in}\mbox{\sl and}\hspace*{0.25in} 
\polyGFib{2m+1}(x)  =  \frac{\lambda^{2m+2}-1}{\lambda^{m}(\lambda^{2}-1)}.\]
{\sl The sign-alternating Lucas polynomials can be written thusly:} 
\[\polyGLucas{2m}(x)  =  \frac{\lambda^{2m}+1}{\lambda^{m}} \hspace*{0.25in}\mbox{\sl and}\hspace*{0.25in}
\polyGLucas{2m+1}(x)  =  \frac{\lambda^{2m+1}+1}{\lambda^{m}(\lambda + 1)}.\]

One can use the preceding forms to confirm, and in fact derive, explicit expressions for the roots of the sign-alternating Fibonacci/Lucas polynomials, cf.\  \ExplicitCorollaryTwo. 
It appears that these results were first obtained in \cite{HB}. 
These explicit expressions are useful for interpreting a key result of the next section (\GameCorollary). 

\noindent 
{\bf \ExplicitCorollaryTwo}\ \ {\sl The sign-alternating Fibonacci polynomials} $\polyGFib{0}(x)$ {\sl and} $\polyGFib{1}(x)$ {\sl and the sign-alternating Lucas polynomials}  $\polyGLucas{0}(x)$ {\sl and} $\polyGLucas{1}(x)$ {\sl are positive constants and therefore have no roots.  
Now let $k$ be an integer with $k \geq 2$. 
The $\lfloor k/2 \rfloor$ distinct roots of the sign-alternating Fibonacci polynomial} $\polyGFib{k}(x)$ {\sl comprise the set}
\[\left\{4\cos^{2}\left(\frac{j\pi}{k+1}\right)\, \rule[-4.5mm]{0.25mm}{11mm}\ \mbox{$j$ is an integer satisfying } 1 \leq j \leq \lfloor k/2 \rfloor\right\}.\]
{\sl Now suppose $d$ and $r$ are integers with $d$ odd and $r$ nonnegative such that $k=2^{r}d$. 
The $\lfloor k/2 \rfloor$ distinct roots of the sign-alternating Lucas polynomial} $\polyGLucas{k}(x)$ {\sl comprise the set}
\[\left\{4\cos^{2}\left(\frac{j\pi}{k} - \frac{\pi}{2^{r+1}}\right)\, \rule[-4.5mm]{0.25mm}{11mm}\ j \in \mbox{$\left\{\frac{d+(2l-1)}{2}\, \rule[-2.5mm]{0.25mm}{7mm}\ \mbox{$l$ is an integer satisfying } 1 \leq l \leq \lfloor k/2 \rfloor\right\}$}\right\}.\]
{\sl When $k$ is odd (i.e.\ when $r=0$), the preceding set of roots can be re-expressed as}
\[\left\{4\sin^{2}\left(\frac{j\pi}{k}\right)\, \rule[-4.5mm]{0.25mm}{11mm}\ \mbox{$j$ is an integer satisfying } 1 \leq j \leq \lfloor k/2 \rfloor\right\}.\]

\newcommand{\TwoCitiesGraphWithSmallerLabels}[2]{
\setlength{\unitlength}{0.75in}
\begin{picture}(1.65,0.25)
\put(0.25,0){\begin{picture}(1,0)
            \put(0,0.1){\circle*{0.05}}
            \put(-0.2,0.2){\small $\gamma_{1}$}
            \put(-0.20,-0.15){\small #1}
            \put(1,0.1){\circle*{0.05}}
            \put(1,0.2){\small $\gamma_{2}$}
            \put(1.05,-0.15){\small #2}
            \put(0,0.1){\line(1,0){1}}
            \put(0.2,0.1){\vector(1,0){0.1}}
            \put(0.8,0.1){\vector(-1,0){0.1}}
            \put(0.225,0.2){\scriptsize $p$}
            \put(0.71,0.2){\scriptsize $q$}
            \end{picture}}
\end{picture}}

\newcommand{\TwoCitiesFiveMoves}[2]{
\setlength{\unitlength}{0.75in}
\begin{picture}(1.65,0.25)
\put(0.25,0){\begin{picture}(1,0)
            \put(0,0.1){\circle*{0.05}}
            \put(-0.2,0.2){\small $\gamma_{1}$}
            \put(-0.20,-0.15){\small #1}
            \put(1,0.1){\circle*{0.05}}
            \put(1,0.2){\small $\gamma_{2}$}
            \put(1.05,-0.15){\small #2}
            \put(0,0.1){\line(1,0){1}}
            \put(0.2,0.1){\vector(1,0){0.1}}
            \put(0.8,0.1){\vector(-1,0){0.1}}
            \put(0.16,0.2){\tiny $7/2$}
            \put(0.65,0.2){\tiny $8/7$}
            \end{picture}}
\end{picture}}

\vspace*{0.1in}
\noindent 
{\large \bf \S \NumbersGame. A Gibonacci generalization of the Networked-numbers Game on two-node graphs.}  
Our variation of the so-called Networked-numbers Game (NG) will be played on a simple connected graph $\Gamma$ with two nodes labelled $\gamma_{1}$ and $\gamma_{2}$. 
Real numbers are assigned to the nodes of $\Gamma$ in pairs; for an ordered pair of real numbers $(u,v)$, we assume $u$ is assigned to $\gamma_{1}$ and $v$ to $\gamma_{2}$. 
Fix positive real numbers $p$ and $q$, to be thought of as multipliers when certain node-firing moves are applied. 
In ordinary NG-play, these moves are as follows. 
The $\gamma_{1}${\em -node-firing move} replaces $u$ with $-u$ and replaces $v$ with $pu+v$, which we depict as follows: 
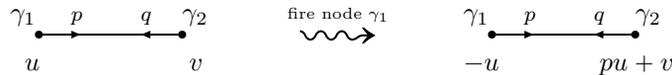
\begin{center}
\TwoCitiesGraphWithSmallerLabels{\hspace*{0.075in}$u$}{$v$} \hspace*{0.2in}{\huge $\stackrel{\mbox{\tiny fire node $\gamma_{1}$}}{\rightsquigarrow}$}\hspace*{0.2in} \TwoCitiesGraphWithSmallerLabels{$-u$}{\hspace*{-0.225in}$pu+v$}
\end{center}
This firing move can only be applied {\em legally} in NG-play if $u>0$. 
Similarly, the $\gamma_{2}${\em -node-firing move} replaces $u$ with $u+qv$ and replaces $v$ with $-v$, which we depict as follows: 
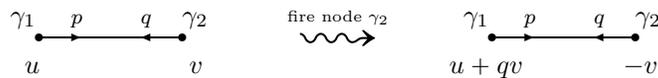
\begin{center}
\TwoCitiesGraphWithSmallerLabels{\hspace*{0.075in}$u$}{$v$} \hspace*{0.2in}{\huge $\stackrel{\mbox{\tiny fire node $\gamma_{2}$}}{\rightsquigarrow}$}\hspace*{0.2in} \TwoCitiesGraphWithSmallerLabels{\hspace*{-0.075in}$u+qv$}{\hspace*{-0.1in}$-v$}
\end{center}
This firing move can only be applied {\em legally} if $v>0$. 

The Networked-numbers Game on two-node graphs begins with a player choosing a pair $(a,b)$ of nonnegative numbers (at least one of which is nonzero) to assign to the nodes of our graph, then choosing $\gamma_{1}$ or $\gamma_{2}$ to apply a legal firing move, and then repeating the preceding step as long as a firing move can be applied legally to the result of the previous firing.
What we call the $(\alpha,\beta)${\em -seeded Gibonacci game} varies the NG only by modifying the initial firing move. 
If the initial move is to fire $\gamma_{1}$, then $(a,b)$ is replaced with $(-\alpha a - q(\alpha-\beta) b\, ,\, p\beta a+ \alpha b)$, but all subsequent node-firings are unmodified and conform to the usual firing rules of NG-play.  
\begin{center}
\TwoCitiesGraphWithSmallerLabels{\hspace*{0.075in}$a$}{$b$} \hspace*{0.2in}{\huge $\stackrel{\mbox{\tiny modified $\gamma_{1}$-firing}}{\rightsquigarrow}$}\hspace*{0.7in} \TwoCitiesGraphWithSmallerLabels{\hspace*{-0.7in}$-\alpha\, a - q(\alpha-\beta)\, b$}{\hspace*{-0.2in}$p\beta\, a+ \alpha\, b$}
\end{center}
Similarly, the modified initial $\gamma_{2}$ firing move replaces $(a,b)$ with $(\alpha a + q\beta b\, ,\,-p(\alpha-\beta) a - \alpha b)$. 
\begin{center}
\hspace*{-0.5in}\TwoCitiesGraphWithSmallerLabels{\hspace*{0.075in}$a$}{$b$} \hspace*{0.2in}{\huge $\stackrel{\mbox{\tiny modified $\gamma_{2}$-firing}}{\rightsquigarrow}$}\hspace*{0.3in} \TwoCitiesGraphWithSmallerLabels{\hspace*{-0.3in}$\alpha\, a + q\beta\, b$}{\hspace*{-0.4in}$-p(\alpha-\beta)\, a - \alpha\, b$}
\end{center}
Note that each of these modified rules reduces to the usual NG rules for initial node-firings when $\alpha = \beta = 1$; more importantly, they ``seed'' certain initial values into game play so that we can utilize sign-alternating Gibonacci polynomials. 

\noindent 
{\bf \GameExample}\ \ In this example, we take $\alpha := 5$, $\beta := 2$, $p := \frac{7}{2}$, and $q := \frac{8}{7}$. 
Viewing $a$ and $b$ as generic positive real numbers, here, is how the Gibonacci game proceeds when we fire $\gamma_{1}$ first:  
\begin{center}
\begin{tabular}{ccc}
\rule[-6.5mm]{0mm}{14mm}\hspace*{-0.0in}\TwoCitiesFiveMoves{\hspace*{0.05in}$a$}{\hspace*{0.0in}$b$} 
& \hspace*{0.0in}{\huge $\stackrel{\mbox{\tiny modified $\gamma_{1}$-firing}}{\rightsquigarrow}$}\hspace*{0.2in} 
& \TwoCitiesFiveMoves{\hspace*{-0.3in}$-5a-\frac{8}{7}\cdot{3b}$}{\hspace*{-0.25in}$\frac{7}{2}\cdot{2a} + 5b$}\\ 
\rule[-6.5mm]{0mm}{14mm} & \hspace*{0.0in}{\huge $\rightsquigarrow$}\hspace*{0.2in} 
& \TwoCitiesFiveMoves{\hspace*{-0.15in}$3a+\frac{16}{7}b$}{\hspace*{-0.2in}$-7a-5b$}\\ 
\rule[-6.5mm]{0mm}{14mm} & \hspace*{0.0in}{\huge $\rightsquigarrow$}\hspace*{0.2in} 
& \TwoCitiesFiveMoves{\hspace*{-0.25in}$-3a-\frac{16}{7}b$}{\hspace*{-0.1in}$\frac{7}{2}a+3b$}\\
\rule[-6.5mm]{0mm}{14mm} & \hspace*{0.0in}{\huge $\rightsquigarrow$}\hspace*{0.2in} 
& \TwoCitiesFiveMoves{\hspace*{-0.075in}$a+\frac{8}{7}b$}{\hspace*{-0.2in}$-\frac{7}{2}a-3b$}\\ 
\rule[-6.5mm]{0mm}{14mm} & \hspace*{0.0in}{\huge $\rightsquigarrow$}\hspace*{0.2in} 
& \TwoCitiesFiveMoves{\hspace*{-0.2in}$-a-\frac{8}{7}b$}{\hspace*{-0.05in}$b$}\\
\rule[-6.5mm]{0mm}{14mm} & \hspace*{0.0in}{\huge $\rightsquigarrow$}\hspace*{0.2in} 
& \TwoCitiesFiveMoves{\hspace*{-0.0in}$-a$}{\hspace*{-0.15in}$-b$} 
\end{tabular}
\end{center}
\noindent 
So, this game terminates in six moves independent of any specific choice for $a$ and $b$.  
Here is how the Gibonacci game proceeds when we fire $\gamma_{2}$ first: 
\begin{center}
\begin{tabular}{ccc}
\rule[-6.5mm]{0mm}{14mm}\hspace*{0.0in}\TwoCitiesFiveMoves{\hspace*{0.05in}$a$}{\hspace*{0.0in}$b$} 
& \hspace*{0.0in}{\huge $\stackrel{\mbox{\tiny modified $\gamma_{2}$-firing}}{\rightsquigarrow}$}\hspace*{0.2in} 
& \TwoCitiesFiveMoves{\hspace*{-0.3in}$5a+\frac{8}{7}\cdot{2b}$}{\hspace*{-0.4in}$-\frac{7}{2}\cdot{3a} - 5b$}\\
\rule[-6.5mm]{0mm}{14mm} & \hspace*{0.0in}{\huge $\rightsquigarrow$}\hspace*{0.2in} 
& \TwoCitiesFiveMoves{\hspace*{-0.25in}$-5a-\frac{16}{7}b$}{\hspace*{-0.1in}$7a+3b$}\\ 
\rule[-6.5mm]{0mm}{14mm} & \hspace*{0.0in}{\huge $\rightsquigarrow$}\hspace*{0.2in} 
& \TwoCitiesFiveMoves{\hspace*{-0.125in}$3a+\frac{8}{7}b$}{\hspace*{-0.2in}$-7a-3b$}\\
\rule[-6.5mm]{0mm}{14mm} & \hspace*{0.0in}{\huge $\rightsquigarrow$}\hspace*{0.2in} 
& \TwoCitiesFiveMoves{\hspace*{-0.225in}$-3a-\frac{8}{7}b$}{\hspace*{-0.1in}$\frac{7}{2}a+b$} 
\end{tabular}
\end{center} 

\begin{center}
\begin{tabular}{ccc}
\rule[-6.5mm]{0mm}{14mm} & \hspace*{0.0in}{\huge $\rightsquigarrow$}\hspace*{0.2in} 
& \TwoCitiesFiveMoves{\hspace*{0.1in}$a$}{\hspace*{-0.2in}$-\frac{7}{2}a-b$}\\
\rule[-6.5mm]{0mm}{14mm} & \hspace*{0.0in}{\huge $\rightsquigarrow$}\hspace*{0.2in} 
& \TwoCitiesFiveMoves{\hspace*{-0.0in}$-a$}{\hspace*{-0.15in}$-b$} 
\end{tabular}
\end{center}
As with the $\gamma_{1}$-first game, this game terminates in six moves independent of any specific choice for $a$ and $b$. 
Note that both games have the same terminal numbers.\hfill\QED 

Our main questions about $(\alpha,\beta)$-seeded Gibonacci games concern termination. 
We need some further terminology to set up these questions and our answers. 
Since Gibonacci game play depends on the choices made for $\alpha$, $\beta$, $p$, and $q$, we refer to the pairing of our two-node graph $\Gamma$ together with a positive-real-number four-tuple $(\alpha, \beta, p, q)$ as a {\em Gibonacci game graph} $\mathcal{G} = \mathcal{G}(\alpha,\beta,p,q)$. 
We say a pair of real numbers $(a,b)$ is {\em nonzero} if at least one of $a$ or $b$ is not zero, {\em dominant} if both $a$ and $b$ are nonnegative, and {\em strongly dominant} if $a$ and $b$ are both positive.
If a Gibonacci game fails to terminate from some given initial choice of numbers $(a,b)$, we say the game {\em diverges}.  
Analogizing Eriksson \cite{E}, we say a Gibonacci game graph is {\em strongly convergent} if, for any given nonzero dominant pair $(a,b)$, any two Gibonacci games either diverge or else terminate in the same number of moves. 
Here, then, are our main questions about Gibonacci games.

\vspace*{-0.15in}
\begin{enumerate}
\item[{\sl (1)}] Which Gibonacci game graphs are strongly convergent? 

\vspace*{-0.1in}
\item[{\sl (2)}] For a strongly convergent Gibonacci game graph, what can be said about the terminal numbers for different terminating games played from a given choice of initial numbers $(a,b)$?
\end{enumerate}

\vspace*{-0.15in}
\noindent 
These are answered by the following theorem. 
Recall the following notation from \S \RootSection: For any $k \geq 2$, the degree $d_{k} := \lfloor k/2 \rfloor$ polynomial $\polyGab{k}(x)$ has $d_{k}$ distinct positive real roots $\{\zeta_{k,1},\ldots,\zeta_{k,d_{k}}\}$, indexed from smallest to largest. 
Set $r_{k} := \zeta_{k,d_{k}}$. 
In addition, let $B^{\alpha,\beta} := \left\{\begin{array}{cl} 4 & \mbox{if } \alpha/\beta \leq 2\\ \frac{(\alpha/\beta)^{2}}{\alpha/\beta - 1} & \mbox{if } \alpha/\beta > 2\end{array}\right.$. 

\noindent 
{\bf \GameTheorem}\ \ {\sl We take as given some Gibonacci game graph $\mathcal{G} = \mathcal{G}(\alpha,\beta,p,q)$. 
(1) No games terminate if $pq \geq B^{\alpha,\beta}$.  If $pq \in (0,B^{\alpha,\beta})$, then all games terminate.  
Moreover, the game graph $\mathcal{G}$ is strongly convergent if and only if $pq \in \{r_{k}\}_{k \geq 2} \cup [B^{\alpha,\beta},\infty)$. 
(2) Suppose, for some $k \geq 2$, we have $pq=r_{k}$.  
Then every game played from a nonzero dominant pair $(a,b)$ terminates at $\left(q\polyGab{k+1}(pq)\, b,-p\polyGab{k-1}(pq)\, a\right) = \left(-q\polyGab{k-1}(pq)\, b,p\polyGab{k+1}(pq)\, a\right)$ when $k$ is even and at $\left(-\polyGab{k-1}(pq)\, a,\polyGab{k+1}(pq)\, b\right) = \left(\polyGab{k+1}(pq)\, a,-\polyGab{k-1}(pq)\, b\right)$ when $k$ is odd. 
Moreover, game play requires exactly $k+1$ node-firings if $(a,b)$ is strongly dominant and exactly $k$ node-firings otherwise.}

{\em Proof.} 
For convenience, let $\widehat{g}_{l} := \polyGab{l}(pq)$ for any nonnegative integer $l$, and set $\widehat{g}_{-1} := \alpha-\beta$. 
It is easy to verify by induction that, if we fire $\gamma_{1}$ (respectively, $\gamma_{2}$) first from a generic strongly dominant pair $(a,b)$, then the $(\alpha,\beta)$-seeded two-node Gibonacci game proceeds as in \GameFigureOne\ (resp.\ \GameFigureTwo).  
\begin{figure}[t]
\begin{center}
\begin{tabular}{ccc}
\rule[-6.5mm]{0mm}{14mm}\TwoCitiesGraphWithSmallerLabels{\hspace*{0.075in}$a$}{$b$}\hspace*{-0.2in} & \hspace*{0.2in}{\huge $\stackrel{\mbox{\tiny modified $\gamma_{1}$-firing}}{\rightsquigarrow}$}\hspace*{0.2in} & \TwoCitiesGraphWithSmallerLabels{\hspace*{-0.35in}$-\widehat{g}_{0}a-q\widehat{g}_{-1}b$}{\hspace*{-0.25in}$p\widehat{g}_{1}a + \widehat{g}_{0}b$}\\ 
\rule[-6.5mm]{0mm}{14mm} & \hspace*{0.2in}{\huge $\rightsquigarrow$}\hspace*{0.2in} 
& \TwoCitiesGraphWithSmallerLabels{\hspace*{-0.3in}$\widehat{g}_{2}a+q\widehat{g}_{1}b$}{\hspace*{-0.35in}$-p\widehat{g}_{1}a - \widehat{g}_{0}b$}\\ 
\rule[-6.5mm]{0mm}{14mm} & \hspace*{0.2in}{\huge $\rightsquigarrow$}\hspace*{0.2in} 
& \TwoCitiesGraphWithSmallerLabels{\hspace*{-0.35in}$-\widehat{g}_{2}a-q\widehat{g}_{1}b$}{\hspace*{-0.25in}$p\widehat{g}_{3}a + \widehat{g}_{2}b$}\\
& {\tiny $\bullet$} \hspace*{0.1in} {\tiny $\bullet$} \hspace*{0.1in} {\tiny $\bullet$} & \\
\rule[-6.5mm]{0mm}{14mm} & \hspace*{0.2in}{\huge $\rightsquigarrow$}\hspace*{0.2in} 
& \TwoCitiesGraphWithSmallerLabels{\hspace*{-0.325in}$\widehat{g}_{k-2}a+q\widehat{g}_{k-3}b$}{\hspace*{-0.225in}$-p\widehat{g}_{k-3}a - \widehat{g}_{k-4}b$}\\
\rule[-6.5mm]{0mm}{14mm} & \hspace*{0.2in}{\huge $\rightsquigarrow$}\hspace*{0.2in} 
& \TwoCitiesGraphWithSmallerLabels{\hspace*{-0.375in}$-\widehat{g}_{k-2}a-q\widehat{g}_{k-3}b$}{\hspace*{-0.1in}$p\widehat{g}_{k-1}a + \widehat{g}_{k-2}b$}\\ 
\rule[-6.5mm]{0mm}{14mm} & \hspace*{0.2in}{\huge $\rightsquigarrow$}\hspace*{0.2in} 
& \TwoCitiesGraphWithSmallerLabels{\hspace*{-0.325in}$\widehat{g}_{k}a+q\widehat{g}_{k-1}b$}{\hspace*{-0.2in}$-p\widehat{g}_{k-1}a - \widehat{g}_{k-2}b$}\\
\rule[-6.5mm]{0mm}{14mm} & \hspace*{0.2in}{\huge $\rightsquigarrow$}\hspace*{0.2in} 
& \TwoCitiesGraphWithSmallerLabels{\hspace*{-0.375in}$-\widehat{g}_{k}a-q\widehat{g}_{k-1}b$}{\hspace*{0.05in}$p\widehat{g}_{k+1}a + \widehat{g}_{k}b$} 
 \end{tabular}

\vspace*{0.1in}
{\small \GameFigureOne: Part of a Gibonacci game with $\gamma_{1}$ fired first. (Here, $k$ even.)}

\vspace*{-0.2in}
\end{center}
\end{figure}
Suppose $pq \geq B^{\alpha,\beta}$. Then by \GibRootTheorem, $\widehat{g}_{l} > 0$. 
Consult \GameFigures\ to see that from any nonzero dominant initial pair $(a,b)$, no Gibonacci game terminates. 

Now suppose $pq \in (0,B^{\alpha,\beta})$.  
First, consider the case $0 < pq < r_{2} = \alpha/\beta$. 
That is, $\beta pq - \alpha <0$, i.e.\ $\widehat{g}_{2} < 0$. 
Then, $\widehat{g}_{3} = \beta pq -\alpha - \beta < 0$ also. 
Assume for the moment that $\gamma_{1}$ is fired first, so we know $a>0$.  
Then, by consultation with \GameFigureOne, the Gibonacci game terminates after two firings if and only if $\widehat{g}_{2}a+q\widehat{g}_{1}b \leq 0$, i.e.\  $\frac{b}{a} \leq -\widehat{g}_{2}/(q\widehat{g}_{1}) = (\alpha-\beta pq)/(q \beta)$.  
In the case that $\frac{b}{a} > (\alpha-\beta pq)/(q \beta)$, then the Gibonacci game terminates after three firings, since $p\widehat{g}_{3}a + \widehat{g}_{2}b$ is necessarily negative. 
That is, $\gamma_{1}$-first Gibonacci games terminate in two moves if and only if $\frac{b}{a} \leq (\alpha-\beta pq)/(q \beta)$ and otherwise terminate in three moves.  
Similarly see that $\gamma_{2}$-first Gibonacci games terminate in two moves if and only if $\frac{a}{b} \leq (\alpha-\beta pq)/(p \beta)$ and otherwise terminate in three moves. 

Next, consider the case $r_{j-1} < pq < r_{j}$ for $j > 2$. 
By \GibRootTheorem, we know that $\widehat{g}_{l}>0$ for $l \in \{0,1,\ldots,j-1\}$, $\widehat{g}_{j}<0$, and $\widehat{g}_{j+1}<0$. 
Assume for the moment that $\gamma_{1}$ is fired first, so $a>0$. 
Supposing that $j$ is even, we can take $k=j$ in \GameFigureOne. 
Clearly $p\widehat{g}_{j-1}a + \widehat{g}_{j-2}b > 0$ and $p\widehat{g}_{j+1}a + \widehat{g}_{j}b < 0$. 
So our Gibonacci game terminates in exactly $j$ firing moves if and only if $\widehat{g}_{j}a+q\widehat{g}_{j-1}b \leq 0$, i.e.\ $\frac{b}{a} \leq -\widehat{g}_{j}/(q\widehat{g}_{j-1})$, and in exactly $j+1$ firing moves otherwise. 
Next assume that $j$ is odd and take $k-1=j$ in \GameFigureOne. 
In this case, $\widehat{g}_{j-1}a+q\widehat{g}_{j-2}b > 0$ and $\widehat{g}_{j+1}a+q\widehat{g}_{j}b < 0$. 
Then our Gibonacci game terminates in exactly $j$ firing moves if and only if $p\widehat{g}_{j}a+\widehat{g}_{j-1}b \leq 0$, i.e.\ $\frac{b}{a} \leq -p\widehat{g}_{j}/\widehat{g}_{j-1}$, and in exactly $j+1$ firing moves otherwise. 
When $\gamma_{2}$ is fired first, similar analysis shows that all games require $j$ or $j+1$ firing moves, with some games of each length. 

\begin{figure}[t]
\begin{center}
\begin{tabular}{ccc}
\rule[-6.5mm]{0mm}{14mm}\TwoCitiesGraphWithSmallerLabels{\hspace*{0.075in}$a$}{$b$}\hspace*{-0.2in} & \hspace*{0.2in}{\huge $\stackrel{\mbox{\tiny modified $\gamma_{2}$-firing}}{\rightsquigarrow}$}\hspace*{0.2in} & \TwoCitiesGraphWithSmallerLabels{\hspace*{-0.3in}$\widehat{g}_{0}a+q\widehat{g}_{1}b$}{\hspace*{-0.35in}$-p\widehat{g}_{-1}a - \widehat{g}_{0}b$}\\ 
\rule[-6.5mm]{0mm}{14mm} & \hspace*{0.2in}{\huge $\rightsquigarrow$}\hspace*{0.2in} 
& \TwoCitiesGraphWithSmallerLabels{\hspace*{-0.35in}$-\widehat{g}_{0}a-q\widehat{g}_{1}b$}{\hspace*{-0.15in}$p\widehat{g}_{1}a + \widehat{g}_{2}b$}\\ 
\rule[-6.5mm]{0mm}{14mm} & \hspace*{0.2in}{\huge $\rightsquigarrow$}\hspace*{0.2in} 
& \TwoCitiesGraphWithSmallerLabels{\hspace*{-0.3in}$\widehat{g}_{2}a+q\widehat{g}_{3}b$}{\hspace*{-0.25in}$-p\widehat{g}_{1}a - \widehat{g}_{2}b$}\\
& {\tiny $\bullet$} \hspace*{0.1in} {\tiny $\bullet$} \hspace*{0.1in} {\tiny $\bullet$} & \\
\rule[-6.5mm]{0mm}{14mm} & \hspace*{0.2in}{\huge $\rightsquigarrow$}\hspace*{0.2in} 
& \TwoCitiesGraphWithSmallerLabels{\hspace*{-0.4in}$-\widehat{g}_{k-4}a-q\widehat{g}_{k-3}b$}{\hspace*{-0.025in}$p\widehat{g}_{k-3}a + \widehat{g}_{k-2}b$}\\
\rule[-6.5mm]{0mm}{14mm} & \hspace*{0.2in}{\huge $\rightsquigarrow$}\hspace*{0.2in} 
& \TwoCitiesGraphWithSmallerLabels{\hspace*{-0.35in}$\widehat{g}_{k-2}a+q\widehat{g}_{k-1}b$}{\hspace*{-0.125in}$-p\widehat{g}_{k-3}a - \widehat{g}_{k-2}b$}\\ 
\rule[-6.5mm]{0mm}{14mm} & \hspace*{0.2in}{\huge $\rightsquigarrow$}\hspace*{0.2in} 
& \TwoCitiesGraphWithSmallerLabels{\hspace*{-0.4in}$-\widehat{g}_{k-2}a-q\widehat{g}_{k-1}b$}{\hspace*{0.125in}$p\widehat{g}_{k-1}a + \widehat{g}_{k}b$}\\
\rule[-6.5mm]{0mm}{14mm} & \hspace*{0.2in}{\huge $\rightsquigarrow$}\hspace*{0.2in} 
& \TwoCitiesGraphWithSmallerLabels{\hspace*{-0.35in}$\widehat{g}_{k}a+q\widehat{g}_{k+1}b$}{\hspace*{0.025in}$-p\widehat{g}_{k-1}a - \widehat{g}_{k}b$} 
 \end{tabular}

\vspace*{0.1in}
{\small \GameFigureTwo: Part of a Gibonacci game with $\gamma_{2}$ fired first. (Here, $k$ even.)}

\vspace*{-0.2in}
\end{center}
\end{figure}

Therefore, $\mathcal{G}$ can only be strongly convergent if $pq \in \{r_{k}\}_{k \geq 2} \cup [B^{\alpha,\beta},\infty)$. 
When $pq \in [B^{\alpha,\beta},\infty)$, then all games diverge so $\mathcal{G}$ is, by definition, strongly convergent. 
Now suppose $pq = r_{j}$ for some $j \geq 2$. 
Then $\widehat{g}_{l}>0$ for $l \in \{0,1,\ldots,j-1\}$, $\widehat{g}_{j}=0$, and $\widehat{g}_{j+1}<0$. 
Notice that $\widehat{g}_{j+1} = (pq)^{j\, \mbox{\scriptsize mod 2}}\widehat{g}_{j} - \widehat{g}_{j-1} = -\widehat{g}_{j-1}$.
Suppose for the moment that $j$ is even. 
Assuming  $a$ is positive, then in \GameFigureOne, we can take $k=j$. 
if $b=0$, this game terminates in $j$ firing moves with terminal numbers $(\widehat{g}_{j}a+q\widehat{g}_{j-1}b\, ,\, -p\widehat{g}_{j-1}a - \widehat{g}_{j-2}b) = (-q\widehat{g}_{j-1}b\, ,\, p\widehat{g}_{j+1}a) = (q\widehat{g}_{j+1}b\, ,\, -p\widehat{g}_{j-1}a) = (-\widehat{g}_{j}a-q\widehat{g}_{j-1}b\, ,\, p\widehat{g}_{j+1}a + \widehat{g}_{j}b)$. 
Now suppose $b$ is positive. 
Then our game terminates in $j+1$ firing moves with terminal numbers $(-\widehat{g}_{j}a-q\widehat{g}_{j-1}b\, ,\, p\widehat{g}_{j+1}a + \widehat{g}_{j}b) = (-q\widehat{g}_{j-1}b\, ,\, p\widehat{g}_{j+1}a) = (q\widehat{g}_{j+1}b\, ,\, -p\widehat{g}_{j-1}a)$. 
With $a$ and $b$ both positive, we can also consider \GameFigureTwo\ with $k=j$. 
This $\gamma_{2}$-first game will terminate in $j+1$ firing moves with terminal numbers $(\widehat{g}_{j}a+q\widehat{g}_{j+1}b\, ,\, -p\widehat{g}_{j-1}a - \widehat{g}_{j}b) = (q\widehat{g}_{j+1}b\, ,\, -p\widehat{g}_{j-1}a) = (-q\widehat{g}_{j-1}b\, ,\, p\widehat{g}_{j+1}a)$. 
These pairs agree with the terminal numbers of the $\gamma_{1}$-first game. 
In the case that $b$ is positive and $a=0$, then from \GameFigureTwo\ we see that the game terminates in $j$ moves with terminal numbers $(-\widehat{g}_{j-2}a-q\widehat{g}_{j-1}b\, ,\, p\widehat{g}_{j-1}a + \widehat{g}_{j}b) = (-q\widehat{g}_{j-1}b\, ,\, p\widehat{g}_{j+1}a) = (q\widehat{g}_{j+1}b\, ,\, -p\widehat{g}_{j-1}a) = (\widehat{g}_{j}a+q\widehat{g}_{j+1}b\, ,\, -p\widehat{g}_{j-1}a - \widehat{g}_{j}b)$.

The preceding paragraph confirms that $\mathcal{G}$ is strongly convergent when $pq \in \{r_{k}\}_{k \geq 2} \cup [B^{\alpha,\beta},\infty)$, that all games terminate when $pq \in \{r_{k}\}_{k \geq 2}$, and that the terminal numbers are as claimed in the theorem statement when $pq \in \{r_{k}\}_{k \geq 2}$.\hfill\QED

The next result specializes the preceding theorem to the Fibonacci ($\alpha=1$, $\beta=1$) and Lucas ($\alpha=2$, $\beta=1$) cases, with the aid of \ExplicitCorollaryTwo. 
Set $S^{1,1} := \left\{ 4\cos^{2}\left(\frac{\pi}{k+1}\right) \right\}_{k \geq 2}$ and $S^{2,1} := \left\{ 4\cos^{2}\left(\rule[-2mm]{-0.1mm}{4.5mm}\frac{\pi}{2k}\right) \right\}_{k \geq 2}$. 

\noindent
{\bf \GameCorollary}\ \ 
{\sl We take as given some Gibonacci game graph $\mathcal{G} = \mathcal{G}(\alpha,\beta,p,q)$.  
For the Fibonacci (respectively, Lucas) case, we take $\alpha=1$, $\beta=1$ (resp.\ $\alpha=2$, $\beta=1$). 
Then no games terminate if $pq \geq 4$.  If $pq \in (0,4)$, then all games terminate.  
Moreover, the game graph $\mathcal{G}$ is strongly convergent if and only if $pq \in S^{\alpha,\beta} \cup [4,\infty)$. 
Suppose, for some $k \geq 2$, we have $pq=4\cos^{2}\left(\frac{\pi}{k+1}\right)$ (resp. $4\cos^{2}\left(\rule[-2mm]{-0.1mm}{4.5mm}\frac{\pi}{2k}\right)$).  
Then every game played from a nonzero dominant pair $(a,b)$ terminates at $\left(q\polyGab{k+1}(pq)\, b,-p\polyGab{k-1}(pq)\, a\right) = \left(-q\polyGab{k-1}(pq)\, b,p\polyGab{k+1}(pq)\, a\right)$ when $k$ is even and at $\left(-\polyGab{k-1}(pq)\, a,\polyGab{k+1}(pq)\, b\right) = \left(\polyGab{k+1}(pq)\, a,-\polyGab{k-1}(pq)\, b\right)$ when $k$ is odd. 
Moreover, game play requires exactly $k+1$ node-firings if $(a,b)$ is strongly dominant and exactly $k$ node-firings otherwise.}

\vspace*{0.1in}
\noindent 
{\large \bf \S \LatticeSection. Symmetric Gibonaccian ranked posets.}  
We now produce some finite ranked posets 
that generalize the `symmetric Fibonaccian lattices' of \cite{DDMN}. 
The symmetric Fibonaccian lattices have the following salutary properties: (1) They are enumerated by a particular specialization of the sign-alternating Fibonacci polynomials; (2) They have rank generating functions whose coefficients are nicely described by some (mostly) new recursively defined symmetric triangular integer arrays;  and (3) They are naturally related to certain representations of the special linear Lie algebras. 
Here, we demonstrate that properties (1) and (2) generalize to the symmetric Gibonaccian ranked posets introduced below. 

We begin by fixing positive integers $\alpha$, $n$, and $k$. 
For reasons that will be explained shortly, we require that $\beta=1$. 
Declare that 
{\small \[R^{\mbox{\tiny $\alpha$-Gib}}(n,k) := \left\{T=(T_{1},\ldots,T_{k})\, \rule[-8mm]{0.2mm}{17mm}\, \begin{array}{c}\CircleInteger{1}T_{j} \in \{\mbox{\small $(j-1)n+1,(j-1)n+2,\ldots,jn$}\}\ \mbox{for all}\ j \in \{1,\ldots,k\},\\ \CircleInteger{2}T_{j+1} \ne T_{j}+1\ \mbox{for all}\ j \in \{1,\ldots,k-1\},\\
\mbox{and}\ \CircleInteger{3}(T_{1},T_{k}) \not\in \{(1,nk),(2,nk-1),\ldots,(\alpha-1,nk-(\alpha-2))\}\end{array}\right\},\]
} 
a set of positive integer $k$-tuples satisfying certain conditions. 
We refer to the objects of this collection as $\alpha${\em -Gibonaccian strings}.
The conditions\, $\CircleInteger{1}\!$ from the above definition are to be called the {\em coordinate requirements}, the conditions\, $\CircleInteger{2}\!$ are the {\em Fibonacci requirements}, and the conditions\, $\CircleInteger{3}\!$ are the $\alpha${\em -requirements}. 
We are most interested in a certain partial ordering of $\alpha$-Gibonaccian strings. 
But first, we consider their enumeration as an unordered collection via inclusion-exclusion, as the latter method makes a direct connection with the sign-alternating Gibonacci polynomials. 
We state this result here as a challenge for the reader and obtain the result later by different means. 

\noindent 
{\bf \InclusionExclusion}\ \  Use the enumerative method of inclusion-exclusion to demonstrate the following equality: $\left|\rule[-1.5mm]{0mm}{4.75mm}R^{\mbox{\tiny $\alpha$-Gib}}(n,k)\right| = n^{k\, \mbox{\scriptsize mod 2}}\, \overline{G}^{\alpha,1}_{k}(n^{2})$.

Now order the $\alpha$-Gibonaccian strings of $R^{\mbox{\tiny $\alpha$-Gib}}(n,k)$ by {\em reverse component-wise comparison}, i.e.\ we have $S \leq T$ in $R^{\mbox{\tiny $\alpha$-Gib}}(n,k)$ for $S = (S_{1},\ldots,S_{k})$ and $T=(T_{1},\ldots,T_{k})$ if and only if $S_{j} \geq T_{j}$ for any $j \in \{1,\ldots,k\}$. 
Observe that $T$ covers $S$ in the resulting order diagram if and only if there exists some $l \in \{1,\ldots,k\}$ such that $S_{l} = T_{l}+1$ while $S_{j} = T_{j}$ for all $j \ne l$. 
One can see that $R^{\mbox{\tiny $\alpha$-Gib}}(n,k)$ is a connected and self-dual ranked poset. 
We call $R^{\mbox{\tiny $\alpha$-Gib}}(n,k)$ a {\em symmetric $\alpha$-Gibonaccian (ranked) poset} (or {\em SGP} for short). 
Of course, $R^{\mbox{\tiny $\alpha$-Gib}}(n,1)$ is an $n$-element chain. 
For convenience later on, we regard $R^{\mbox{\tiny $\alpha$-Gib}}(n,0)$ to be an $\alpha$-element anti-chain. 
In \RalphaGibFigure, we depict $R^{\mbox{\tiny $3$-Gib}}(4,3)$. 

The next result says that for fixed $\alpha$, SGPs are distributive lattices for all $k$ only when $\alpha$ is $1$ or $2$, which we refer to respectively as the {\em symmetric Fibonaccian lattices} ($\alpha =1$) and {\em symmetric Lucasian lattices} ($\alpha = 2$). 

\noindent
{\bf \RankedProp}\ \ {\sl Assume $n>\alpha$.  Then the connected and self-dual ranked poset} $R^{\mbox{\tiny $\alpha$-Gib}}(n,k)$ {\sl is a distributive lattice for some $k > 1$ if and only if $\alpha \in \{1,2\}$ if and only if} $R^{\mbox{\tiny $\alpha$-Gib}}(n,k)$ {\sl has a unique maximal element for some $k > 1$ if and only if} $R^{\mbox{\tiny $\alpha$-Gib}}(n,k)$ {\sl is a distributive lattice for all $k > 1$.} 

\begin{figure}[tb]
\begin{center}
\RalphaGibFigure: The 48-element symmetric $3$-Gibonaccian poset $R^{\mbox{\tiny $3$-Gib}}(4,3)$.\\
$\RGF(R^{\mbox{\tiny $3$-Gib}}(4,3);q) = 1 + 3q + 6q^{2} + 6q^{3} + 8q^{4} + 8q^{5} + 6q^{6} + 6q^{7} + 3q^{8} + q^{9}$. 

\vspace*{-0.1in}
\setlength{\unitlength}{1.25cm}
\begin{picture}(10,9.5)
%
\put(0,4){\VertexBlankFib}
\put(-0.8,4){\tiny $(4,7,9)$}
\put(0,5){\VertexBlankFib}
\put(-0.8,5){\tiny $(4,6,9)$}
\put(1,5){\VertexBlankFib}
\put(0.225,5){\tiny $(3,7,9)$}
\put(1,6){\VertexBlankFib}
\put(0.225,6){\tiny $(3,6,9)$}
\put(1,7){\VertexBlankFib}
\put(0.225,7){\tiny $(3,5,9)$}
\put(2,6){\VertexBlankFib}
\put(1.225,6){\tiny $(2,7,9)$}
\put(2,7){\VertexBlankFib}
\put(1.225,7){\tiny $(2,6,9)$}
\put(2,8){\VertexBlankFib}
\put(1.225,8){\tiny $(2,5,9)$}
\put(3,7){\VertexBlankFib}
\put(3.1,7.05){\tiny $(1,7,9)$}
\put(3,8){\VertexBlankFib}
\put(3.1,8.05){\tiny $(1,6,9)$}
\put(3,9){\VertexBlankFib}
\put(3.1,9.05){\tiny $(1,5,9)$}
\put(2.5,4){\VertexBlankFib}
\put(1.6,3.8){\tiny $(4,6,10)$}
\put(2.5,3){\VertexBlankFib}
\put(1.6,2.8){\tiny $(4,7,10)$}
\put(2.5,2){\VertexBlankFib}
\put(1.6,2){\tiny $(4,8,10)$}
\put(3.5,6){\VertexBlankFib}
\put(2.6,5.925){\tiny $(3,5,10)$}
\put(3.5,5){\VertexBlankFib}
\put(2.6,4.925){\tiny $(3,6,10)$}
\put(3.5,4){\VertexBlankFib}
\put(3.65,4){\tiny $(3,7,10)$}
\put(3.5,3){\VertexBlankFib}
\put(3.65,3){\tiny $(3,8,10)$}
\put(4.5,7){\VertexBlankFib}
\put(4.6,6.9){\tiny $(2,5,10)$}
\put(4.5,6){\VertexBlankFib}
\put(4.6,5.9){\tiny $(2,6,10)$}
\put(4.5,5){\VertexBlankFib}
\put(4.6,4.9){\tiny $(2,7,10)$}
\put(4.5,4){\VertexBlankFib}
\put(4.6,3.9){\tiny $(2,8,10)$}
\put(5.5,8){\VertexBlankFib}
\put(5.6,8.05){\tiny $(1,5,10)$}
\put(5.5,7){\VertexBlankFib}
\put(5.6,7.05){\tiny $(1,6,10)$}
\put(5.5,6){\VertexBlankFib}
\put(5.6,6.05){\tiny $(1,7,10)$}
\put(5.5,5){\VertexBlankFib}
\put(5.65,5.45){\tiny $(1,8,10)$}
\put(5.7,5.4){\vector(-1,-2){0.175}}
\put(5,3){\VertexBlankFib}
\put(4.15,2.8){\tiny $(4,6,11)$}
\put(5,2){\VertexBlankFib}
\put(4.15,1.8){\tiny $(4,7,11)$}
\put(5,1){\VertexBlankFib}
\put(4.15,0.8){\tiny $(4,8,11)$}
\put(6,5){\VertexBlankFib}
\put(6.1,5.05){\tiny $(3,5,11)$}
\put(6,4){\VertexBlankFib}
\put(6.1,4.05){\tiny $(3,6,11)$}
\put(6,3){\VertexBlankFib}
\put(6.1,3.05){\tiny $(3,7,11)$}
\put(6,2){\VertexBlankFib}
\put(6.1,2.05){\tiny $(3,8,11)$}
\put(8,7){\VertexBlankFib}
\put(8.1,7.05){\tiny $(1,5,11)$}
\put(8,6){\VertexBlankFib}
\put(8.1,6.05){\tiny $(1,6,11)$}
\put(8,5){\VertexBlankFib}
\put(8.1,5.05){\tiny $(1,7,11)$}
\put(8,4){\VertexBlankFib}
\put(7.5,3.75){\tiny $(1,8,11)$}
\put(7.5,2){\VertexBlankFib}
\put(6.65,1.8){\tiny $(4,6,12)$}
\put(7.5,1){\VertexBlankFib}
\put(6.65,0.8){\tiny $(4,7,12)$}
\put(7.5,0){\VertexBlankFib}
\put(6.65,-0.2){\tiny $(4,8,12)$}
\put(8.5,4){\VertexBlankFib}
\put(8.6,3.95){\tiny $(3,5,12)$}
\put(8.5,3){\VertexBlankFib}
\put(8.6,2.95){\tiny $(3,6,12)$}
\put(8.5,2){\VertexBlankFib}
\put(8.6,1.95){\tiny $(3,7,12)$}
\put(8.5,1){\VertexBlankFib}
\put(8.6,0.95){\tiny $(3,8,12)$}
\put(9.5,5){\VertexBlankFib}
\put(9.6,5.05){\tiny $(2,5,12)$}
\put(9.5,4){\VertexBlankFib}
\put(9.6,4.05){\tiny $(2,6,12)$}
\put(9.5,3){\VertexBlankFib}
\put(9.6,3.05){\tiny $(2,7,12)$}
\put(9.5,2){\VertexBlankFib}
\put(9.6,2.05){\tiny $(2,8,12)$}
%
\put(0,0){\qbezier(0,4)(1.5,5.5)(3,7)}
\put(0,0){\qbezier(0,5)(1.5,6.5)(3,8)}
\put(0,0){\qbezier(1,7)(2,8)(3,9)}
\put(0,0){\qbezier(0,4)(0,4.5)(0,5)}
\put(0,0){\qbezier(1,5)(1,6)(1,7)}
\put(0,0){\qbezier(2,6)(2,7)(2,8)}
\put(0,0){\qbezier(3,7)(3,8)(3,9)}
\put(0,0){\qbezier(2.5,2)(2.5,3)(2.5,4)}
\put(0,0){\qbezier(3.5,3)(3.5,4.5)(3.5,6)}
\put(0,0){\qbezier(4.5,4)(4.5,5.5)(4.5,7)}
\put(0,0){\qbezier(5.5,5)(5.5,6.5)(5.5,8)}
\put(0,0){\qbezier(2.5,2)(4,3.5)(5.5,5)}
\put(0,0){\qbezier(2.5,3)(4,4.5)(5.5,6)}
\put(0,0){\qbezier(2.5,4)(4,5.5)(5.5,7)}
\put(0,0){\qbezier(3.5,6)(4.5,7)(5.5,8)}
\put(0,0){\qbezier(5,1)(5,2)(5,3)}
\put(0,0){\qbezier(6,2)(6,3.5)(6,5)}
\put(0,0){\qbezier(8,4)(8,5.5)(8,7)}
\put(0,0){\qbezier(5,1)(5.5,1.5)(6,2)}
\put(0,0){\qbezier(5,2)(5.5,2.5)(6,3)}
\put(0,0){\qbezier(5,3)(5.5,3.5)(6,4)}
\put(0,0){\qbezier(7.5,0)(7.5,1)(7.5,2)}
\put(0,0){\qbezier(8.5,1)(8.5,2.5)(8.5,4)}
\put(0,0){\qbezier(9.5,2)(9.5,3.5)(9.5,5)}
\put(0,0){\qbezier(7.5,0)(8.5,1)(9.5,2)}
\put(0,0){\qbezier(7.5,1)(8.5,2)(9.5,3)}
\put(0,0){\qbezier(7.5,2)(8.5,3)(9.5,4)}
\put(0,0){\qbezier(8.5,4)(9,4.5)(9.5,5)}
\put(0,0){\qbezier(3,9)(5.5,8)(8,7)}
\put(0,0){\qbezier(3,8)(5.5,7)(8,6)}
\put(0,0){\qbezier(3,7)(5.5,6)(8,5)}
\put(0,0){\qbezier(5.5,5)(6.75,4.5)(8,4)}
\put(0,0){\qbezier(2,8)(3.25,7.5)(4.5,7)}
\put(0,0){\qbezier(2,7)(3.25,6.5)(4.5,6)}
\put(0,0){\qbezier(2,6)(3.25,5.5)(4.5,5)}
\put(0,0){\qbezier(1,7)(4.75,5.5)(8.5,4)}
\put(0,0){\qbezier(1,6)(4.75,4.5)(8.5,3)}
\put(0,0){\qbezier(1,5)(4.75,3.5)(8.5,2)}
\put(0,0){\qbezier(3.5,3)(6,2)(8.5,1)}
\put(0,0){\qbezier(0,5)(3.75,3.5)(7.5,2)}
\put(0,0){\qbezier(0,4)(3.75,2.5)(7.5,1)}
\put(0,0){\qbezier(2.5,2)(5,1)(7.5,0)}
\end{picture}

\vspace*{-0.2in}
\end{center}
\end{figure}

{\em Proof.} Throughout, we assume $k \geq 2$.  
When $\alpha > 2$, the string $(2,n+1,2n+1,\ldots,n(k-1)+1,nk)$ is in $R^{\mbox{\tiny $\alpha$-Gib}}(n,k)$ but neither $(1,n+1,2n+1,\ldots,n(k-1)+1,nk)$ nor $(2,n+1,2n+1,\ldots,n(k-1)+1,nk-1)$ is. 
So, $R^{\mbox{\tiny $\alpha$-Gib}}(n,k)$ has at least two maximal elements, namely $(1,n+1,\ldots,n(k-1)+1)$ and $(2,n+1,2n+1,\ldots,n(k-1)+1,nk)$, and cannot be a distributive lattice. 

It remains to argue that for any $\alpha \in \{1,2\}$ and any $k \geq 2$, the poset $R^{\mbox{\tiny $\alpha$-Gib}}(n,k)$ is a distributive lattice (in which case $R^{\mbox{\tiny $\alpha$-Gib}}(n,k)$ will have a unique maximal element). 
By Proposition 10.2 of \cite{DonDiamond}, it suffices to show that $R^{\mbox{\tiny $\alpha$-Gib}}(n,k)$ is closed under reverse-componentwise joins and meets. 
Say $S = (S_{1},\ldots,S_{k})$ and $T=(T_{1},\ldots,T_{k})$ are strings in $R^{\mbox{\tiny $\alpha$-Gib}}(n,k)$. 
First, we argue that $\left(\rule[-1.5mm]{-0.8mm}{4mm}\min(S_{i},T_{i})\right)_{i\in\{1,\ldots,k\}}$ is also in $R^{\mbox{\tiny $\alpha$-Gib}}(n,k)$. 
As each of $S$ and $T$ meet all coordinate requirements for membership in $R^{\mbox{\tiny $\alpha$-Gib}}(n,k)$, then so does $\left(\rule[-1.5mm]{0mm}{4mm}\min(S_{i},T_{i})\right)_{i\in\{1,\ldots,k\}}$. 
Now we check the Fibonacci requirements. 
Well, if $\min(S_{i},T_{i}) + 1 = \min(S_{i+1},T_{i+1})$, then $\min(S_{i},T_{i}) = ni$ and $\min(S_{i+1},T_{i+1}) = ni+1$. 
Since $ni$ is the largest number allowed in the $i^{\mbox{\tiny th}}$ coordinate of any element of $R^{\mbox{\tiny $\alpha$-Gib}}(n,k)$, then the only way we get $\min(S_{i},T_{i}) = ni$ is if $S_{i} = ni = T_{i}$. 
Since one of $S_{i+1}$ or $T_{i+1}$ must be $ni+1$, then one of $S$ or $T$ is not in $R^{\mbox{\tiny $\alpha$-Gib}}(n,k)$, contrary to our hypothesis. 
That is, $\left(\rule[-1.5mm]{-0.8mm}{4mm}\min(S_{i},T_{i})\right)_{i\in\{1,\ldots,k\}}$ meets all Fibonacci requirements. 
Last, we check that $\left(\rule[-1.5mm]{-0.8mm}{4mm}\min(S_{i},T_{i})\right)_{i\in\{1,\ldots,k\}}$ meets all $\alpha$-requirements. 
When $\alpha=1$, these requirements are empty. 
So, consider $\alpha=2$. 
We must show that if $\min(S_{1},T_{1}) = 1$, then $\min(S_{k},T_{k}) < nk$. 
Well, if $S_{1} = 1$, then $S_{k} < nk$, and therefore $\min(S_{k},T_{k}) < nk$, and similarly if $T_{1} = 1$. 
That is, $\left(\rule[-1.5mm]{-0.8mm}{4mm}\min(S_{i},T_{i})\right)_{i\in\{1,\ldots,k\}}$ meets all $\alpha$-requirements. 
Second, we must argue that $\left(\rule[-1.5mm]{-0.8mm}{4mm}\max(S_{i},T_{i})\right)_{i\in\{1,\ldots,k\}}$ meets all coordinate, Fibonacci, and $\alpha$-requirements. 
This follows by reasoning entirely similar to the $\left(\rule[-1.5mm]{-0.8mm}{4mm}\min(S_{i},T_{i})\right)_{i\in\{1,\ldots,k\}}$ case.\hfill\QED

We remark that when $n=3$ and $\alpha=1$, the sequence of symmetric Fibonaccian lattice sizes, starting with $k=0$, is $1, 3, 8, 21, 55, 144, \ldots$, coinciding with the Fibonacci subsequence $\{f_{2m+1}\}_{m \geq 0}$ (cf.\ OEIS A001906 \cite{OEIS}).  
We call the numbers of this latter subsequence the {\em symmetric Fibonacci numbers.} 
When $n=3$ and $\alpha=2$, it follows from \EnumerativeTheorem\ below that the analogous sequence of symmetric Lucasian lattice sizes is $2, 3, 7, 18, 47, 123, \ldots$, coinciding with the Lucas subsequence $\{\myl_{2m}\}_{m \geq 0}$ (cf.\ OEIS A005248 \cite{OEIS}). 
We call the numbers of this latter subsequence the {\em symmetric Lucas numbers}. 

Next, we consider a family of triangular arrays indexed by pairs of integers $(\alpha;n)$. 
Our interest was inspired by the ``$(1;3)$'' array presented and investigated in \cite{KK}. 
The $n^{\mbox{\tiny \em th}}$ {\em symmetric $(\alpha;n)$-Gibonacci triangle} $\mymathfrak{A}^{(\alpha;n)} = (\mya^{(\alpha;n)}_{k,r})$ is defined recursively as follows. 
For each nonnegative integer $k$, set $\mathcal{I}_{n,k} := \{-k(n-1),-k(n-1)+2,\ldots,k(n-1)-2,k(n-1)\}$, which is to be thought of as an indexing set for row $k$, and for integers $r \not\in \mathcal{I}_{n,k}$, set $\mya^{(\alpha;n)}_{k,r} := 0$. 
For initial array values, take $\mya^{(\alpha;n)}_{0,0} := \alpha$ with $\mya^{(\alpha;n)}_{1,r} := 1$ for each $r \in \mathcal{I}_{n,1}$. 
Then for $k \geq 2$ and $r \in \mathbb{Z}$, let $\displaystyle \mya^{(\alpha;n)}_{k,r} := \left(\sum_{s \in \mathcal{I}_{n,1}}\mya^{(\alpha;n)}_{k-1,r+s}\right) - \mya^{(\alpha;n)}_{k-2,r}$. 
When $r \in \mathcal{I}_{n,k}$ for a nonnegative integer $k$, we say the array entry $\mya^{(\alpha;n)}_{k,r}$ is {\em regular}. 
Here, for example, is part of $\mymathfrak{A}^{(3;4)}$:  
{
\begin{center}
\renewcommand{\baselinestretch}{1.2}
\footnotesize
\begin{tabular}{ccccccccccccccccccc}
 & & & & & & & & & 3 & & & & & & & & & \\
 & & & & & & 1 & & 1 & & 1 & & 1 & & & & & & \\
 & & & 1 & & 2 & & 3 & & 1 & & 3 & & 2 & & 1 & & & \\
1 & & 3 & & 6 & & 6 & & 8 & & 8 & & 6 & & 6 & & 3 & & 1 
\end{tabular}

\vspace*{0.1in}
ETC.
\end{center}
}
\noindent
We define a polynomial $\displaystyle {A}^{(\alpha;n)}_{k}(x)$ by $\displaystyle {A}^{(\alpha;n)}_{k}(x) := \sum_{r \in \mathcal{I}_{n,k}}\mya^{(\alpha;n)}_{k,r}x^{\frac{1}{2}\left(\rule[-1.25mm]{0mm}{3.5mm}k(n-1)-r\right)}$, and set $\displaystyle {A}^{(\alpha;n)}_{-1}(x) := 0$. 
So, $\displaystyle {A}^{(\alpha;n)}_{k}(1)$ is the sum of the regular entries of the $k^{\mbox{\tiny th}}$ row of $\mymathfrak{A}^{(\alpha;n)}$. 

We now explain why we take $\beta=1$.  
For $k \geq 2$ the first (resp.\ last) regular entry on the $k^{\mbox{\tiny th}}$ row of $\mymathfrak{A}^{(\alpha;n)}$ agrees with the first (resp.\ last) regular entry of  the preceding row. 
Since we declare each regular entry of the $1^{\mbox{\tiny st}}$ row to be unity, then the first and last regular entries of all subsequent rows are also unity. 
The obvious way to modify $\mymathfrak{A}^{(\alpha;n)}$ in order to account for a $\beta$ value larger than one is to replace all $1^{\mbox{\tiny st}}$ row regular entries with $\beta$. 
Then, all subsequent rows would begin and end with $\beta$. 
However, our aim is to generalize symmetric Fibonaccian distributive lattices, which have exactly one element of maximal and one of minimal rank. 
So, here we will require that $\beta=1$. 

The following proposition further helps justify some of our language and conventions. 

\noindent 
{\bf \AlphaNProp}\ \ {\sl For each $k \geq 0$ and $r \in \mathcal{I}_{n,k}$, we have} $\mya^{(\alpha;n)}_{k,r} = \mya^{(\alpha;n)}_{k,-r}${\sl , so the symmetric $(\alpha;n)$-Gibonacci triangle} $\mymathfrak{A}^{(\alpha;n)}$ {\sl indeed has symmetric rows. 
Moreover, all regular entries of} $\mymathfrak{A}^{(\alpha;n)}$ {\sl are positive integers if and only if $n > \alpha$.}

{\em Proof.} The first claim of the proposition statement follows easily by induction on the row numbers $k$.  
For the second claim, that $n > \alpha$ is necessary for positivity of all regular entries of $\mymathfrak{A}^{(\alpha;n)}$ follows from the simple observation that $\mya^{(\alpha;n)}_{2,0} = n-\alpha$. 

Next, we show by induction on row numbers $k$ that $n > \alpha$ is also sufficient. 
Of course, all $n$ of the regular entries on row 1 are unity and therefore positive. 
All regular entries on row 2 are equal to $n$ except that $\mya^{(\alpha;n)}_{2,0} = n-\alpha$, so all of these entries are positive under the hypothesis that $n > \alpha$. 
Now suppose that for some $k > 2$, we know that all regular entries on row $k'$ are positive, if $1 \leq k' < k$.  
Our aim is to show that a generic $k^{\mbox{\tiny th}}$ row regular entry $\mya^{(\alpha;n)}_{k,r}$ is positive. 
By symmetry of the array $\mymathfrak{A}^{(\alpha;n)}$, we can, without loss of generality, assume $r \geq 0$. 
In the formula  
$\displaystyle \mya^{(\alpha;n)}_{k,r} := \left(\sum_{s \in \mathcal{I}_{n,1}}\mya^{(\alpha;n)}_{k-1,r+s}\right) - \mya^{(\alpha;n)}_{k-2,r}$,  
all of the summands of the form $\mya^{(\alpha;n)}_{k-1,r+s}$ are nonnegative (by our inductive hypothesis) and at least one of them is positive. 
Now, either $\mya^{(\alpha;n)}_{k-2,r} = 0$ or $\mya^{(\alpha;n)}_{k-2,r} > 0$.  
If $\mya^{(\alpha;n)}_{k-2,r} = 0$, we can conclude that the $k^{\mbox{\tiny th}}$ row regular entry $\mya^{(\alpha;n)}_{k,r}$ is positive, completing the induction argument. 
Now consider the case that $\mya^{(\alpha;n)}_{k-2,r} > 0$.  
It is easy to see that $\mya^{(\alpha;n)}_{k-2,r}>0$ if and only if $r \in \mathcal{I}_{n,k-2}$ if and only if $r+(n-1) \in \mathcal{I}_{n,k-1}$ if and only if $\mya^{(\alpha;n)}_{k-1,r+(n-1)}>0$.  
The appearance of ``$+\mya^{(\alpha;n)}_{k-2,r}$'' in the formula for $\mya^{(\alpha;n)}_{k-1,r+(n-1)}$ cancels the ``$-\mya^{(\alpha;n)}_{k-2,r}$'' appearing in the formula for $\mya^{(\alpha;n)}_{k,r}$, but at the expense of introducing another (potentially) negative term, namely ``$\mya^{(\alpha;n)}_{k-3,r+(n-1)}$.'' 
Either $\mya^{(\alpha;n)}_{k-3,r+(n-1)} = 0$ or $-\mya^{(\alpha;n)}_{k-3,r+(n-1)} > 0$. 
If $\mya^{(\alpha;n)}_{k-3,r+(n-1)} = 0$, then, as before, we can conclude that  $\mya^{(\alpha;n)}_{k,r}>0$, completing the induction argument. 
So, suppose $\mya^{(\alpha;n)}_{k-3,r+(n-1)} > 0$. 
Well, again observe that $\mya^{(\alpha;n)}_{k-3,r+(n-1)} > 0$ if and only if $r+(n-1) \in \mathcal{I}_{n,k-3}$ if and only if $r+2(n-1) \in \mathcal{I}_{n,k-2}$ if and only if $\mya^{(\alpha;n)}_{k-2,r+2(n-1)}>0$. 
The appearance of ``$+\mya^{(\alpha;n)}_{k-3,r+(n-1)}$'' in the formula for $\mya^{(\alpha;n)}_{k-2,r+2(n-1)}$ cancels the ``$-\mya^{(\alpha;n)}_{k-3,r+(n-1)}$'' now appearing in the formula for $\mya^{(\alpha;n)}_{k,r}$, but at the expense of introducing another (potentially) negative term, namely ``$\mya^{(\alpha;n)}_{k-4,r+2(n-1)}$.'' 
To complete the induction argument, repeat this process until some $\mya^{(\alpha;n)}_{k-2i,r+i(n-1)}=0$.\hfill\QED

In view of the preceding result, from here on, we make the simplifying hypothesis that $n > \alpha$. 
The main results of this section are expressed in \EnumerativeTheorem\ and relate the cardinality and rank sizes of $R^{\mbox{\tiny $\alpha$-Gib}}(n,k)$ respectively to the sign-alternating Gibonacci polynomial $\overline{G}^{\alpha,1}_{k}(x)$ and to the symmetric $(\alpha;n)$-Gibonacci triangle $\mymathfrak{A}^{(\alpha;n)}$. 
These results are stated as enumerative identities and as equalities of certain polynomials in the variable $q$. 
We will use the notation $[m] := (q^{m}-1)/(q-1)$ and call $[m]$ a $q$-integer. 

The ranked poset $R^{\mbox{\tiny $\alpha$-Gib}}(n,k)$ has $M = (1,n+1,\ldots,(k-1)n+1)$ as its unique element of maximal rank, $N = (n,2n,\ldots,kn)$ as its unique element of minimal rank, and length $k(n-1)$. 
Then, one can see that the rank function $\rho: R^{\mbox{\tiny $\alpha$-Gib}}(n,k) \longrightarrow \{0,\ldots,k(n-1)\}$ is given by \[\rho(T) = k(n-1) - \sum_{i=1}^{k}(T_{i}-(i-1)n-1) = \frac{1}{2}k(k+1)n - \sum_{i=1}^{k} T_{i}.\]

\noindent
{\bf \RGFTheorem}\ \ {\sl Let $n$ and $k$ be integers, both larger than $1$. Let ${H}^{(\alpha;n)}_{k}(q)$ denote the rank generating function} $\RGF(R^{\mbox{\tiny $\alpha$-Gib}}(n,k);q)$. 
{\sl By convention, ${A}^{(\alpha;n)}_{0}(q) = \alpha = {H}^{(\alpha;n)}_{0}(q)$ and ${A}^{(\alpha;n)}_{1}(q) = [n] = {H}^{(\alpha;n)}_{1}(q)$.} 
{\sl We have (1) ${A}^{(\alpha;n)}_{k}(q) = [n]{A}^{(\alpha;n)}_{k-1}(q)-q^{n-1}{A}^{(\alpha;n)}_{k-2}(q)$.} 
{\sl When $\alpha=1$, we have (2) ${H}^{(1;n)}_{k}(q) = {A}^{(1;n)}_{k}(q)$.}  
{\sl For general $\alpha$, (3) ${H}^{(\alpha;n)}_{k}(q) = [n]{H}^{(1;n)}_{k-1}(q) - \left(\rule[-1.5mm]{0mm}{4mm}[n]-[n-\alpha]\right){H}^{(1;n)}_{k-2}(q)$.} 
{\sl Consequently, we obtain (4) ${H}^{(\alpha;n)}_{k}(q) = [n]{H}^{(\alpha;n)}_{k-1}(q) - q^{n-1}{H}^{(\alpha;n)}_{k-2}(q)$ and (5)} ${H}^{(\alpha;n)}_{k}(q) = {A}^{(\alpha;n)}_{k}(q)$.\\ 
{\bf \CountingTheorem}\ \ {\sl Continuing, set} 
$\mathscr{H}^{\alpha}_{n,k} := \left|\rule[-1.5mm]{0mm}{4.75mm}R^{\mbox{\tiny $\alpha$-Gib}}(n,k)\right| = {H}^{(\alpha;n)}_{k}(1)$, 
$\mathscr{G}^{\alpha}_{n,k} := n^{k\, \mbox{\scriptsize mod 2}}\, \overline{G}^{\alpha,1}_{k}(n^{2})$, 
{\sl and} $\mathscr{A}^{\alpha}_{n,k} := {A}^{(\alpha;n)}_{k}(1)$. 
{\sl Then for each $\mathscr{X} \in \{\mathscr{A},\mathscr{G},\mathscr{H}\}$ we have $\mathscr{X}^{\alpha}_{n,k} = n\mathscr{X}^{\alpha}_{n,k-1} - \mathscr{X}^{\alpha}_{n,k-2}$, with $\mathscr{X}^{\alpha}_{n,0} = \alpha$ and $\mathscr{X}^{\alpha}_{n,1} = n$. 
When $n=2$, then, under our simplifying hypothesis, necessarily $\alpha=1$, and we have $\mathscr{X}^{1}_{2,k} = k+1$; when $n > 2$, we have} 
$\displaystyle \mathscr{X}^{\alpha}_{n,k} = \frac{\ \myr_{2}^{k}(n-\alpha \myr_{1}) - \myr_{1}^{k}(n-\alpha \myr_{2})\ }{\myr_{2}-\myr_{1}}$,   
{\sl where $\myr_{2}$ is the largest and $\myr_{1}$ is the smallest of the two distinct real roots of $x^{2}-nx+1$.}

{\em Proof.}  The identity of {\bf A.}{\sl (1)} is routine and follows by applying the defining recurrence for the symmetric Gibonacci triangle $\mymathfrak{A}^{(\alpha;n)}$ to the definition of ${A}^{(\alpha;n)}_{k}(q)$.  
The identity ${H}^{(1;n)}_{k}(q) = {A}^{(1;n)}_{k}(q)$ of {\bf A.}{\sl (2)} follows from Theorem 5.1 of \cite{DDMN}. 
For  {\bf A.}{\sl (3)}, consider that an $\alpha$-Gibonacci string $S = (S_{1},\ldots,S_{k})$ in $R^{\mbox{\tiny $\alpha$-Gib}}(n,k)$ has $S_{1} \in \{1,2,\ldots,n\}$. 
If we fix $S_{1} = 1$ (and assume $S_{k} \ne nk$ if $\alpha > 1$), then we can identify $S$ with a string of the form $S' = (S_{2}-n,S_{3}-n,\ldots,S_{k}-n)$ in $R^{\mbox{\tiny $\alpha$-Gib}}(n,k-1)$, and the rank of $S$ is $q^{n-1}$ times the rank of $S'$. 
So if $\alpha=1$, then $\displaystyle \sum_{S\, \mbox{\scriptsize with}\, S_{1}=1}q^{\rho(S)} = q^{n-1}{H}^{(1;n)}_{k-1}(q)$. 
But, if $\alpha > 1$, then we must throw out of said sum all strings with $S_{k} = nk$, hence $\displaystyle \sum_{S\, \mbox{\scriptsize with}\, S_{1}=1}q^{\rho(S)} = q^{n-1}{H}^{(1;n)}_{k-1}(q) - q^{n-1}{H}^{(1;n)}_{k-2}(q)$. 
Now fix $S_{1} = 2$ (and assume $S_{k} \ne nk-1$ if $\alpha > 2$), then we can identify $S$ with a string of the form $S' = (S_{2}-n,S_{3}-n,\ldots,S_{k}-n)$ in $R^{\mbox{\tiny $\alpha$-Gib}}(n,k-1)$, and the rank of $S$ is $q^{n-2}$ times the rank of $S'$. 
So if $\alpha \leq 2$, then $\displaystyle \sum_{S\, \mbox{\scriptsize with}\, S_{1}=2}q^{\rho(S)} = q^{n-2}{H}^{(1;n)}_{k-1}(q)$. 
But, if $\alpha > 2$, then we must throw out of said sum all strings with $S_{k} = nk-1$, hence $\displaystyle \sum_{S\, \mbox{\scriptsize with}\, S_{1}=2}q^{\rho(S)} = q^{n-2}{H}^{(1;n)}_{k-1}(q) - q^{n-2}{H}^{(1;n)}_{k-2}(q)$. 
Continuing in this way, we see that ${H}^{(\alpha;n)}_{k}(q) = \left(\rule[-1.5mm]{0mm}{4mm}q^{n-1}+q^{n-2}+\cdots+1\right){H}^{(1;n)}_{k-1}(q) - \left(\rule[-1.5mm]{0mm}{4mm}q^{n-1}+\cdots+q^{n-\alpha}\right){H}^{(1;n)}_{k-2}(q) = [n]{H}^{(1;n)}_{k-1}(q) - \left(\rule[-1.5mm]{0mm}{4mm}[n]-[n-\alpha]\right){H}^{(1;n)}_{k-2}(q)$. 
Identity {\bf A.}{\sl (4)} is easily obtained by combining the identities {\bf A.}{\sl (1)}, {\bf A.}{\sl (2)}, and {\bf A.}{\sl (3)}. 
Since each of  ${H}^{(\alpha;n)}_{k}(q)$ and ${A}^{(\alpha;n)}_{k}(q)$ satisfy the same recurrence with the same initial conditions by parts {\bf A.}{\sl (1)} and {\bf A.}{\sl (4)}, we get ${H}^{(\alpha;n)}_{k}(q) = {A}^{(\alpha;n)}_{k}(q)$.

For {\bf B}, it is routine to verify that the $\mathscr{G}^{\alpha}_{n,k}$'s satisfy the claimed recurrence. 
That $\mathscr{A}^{\alpha}_{n,k} = n\mathscr{A}^{\alpha}_{n,k-1} - \mathscr{A}^{\alpha}_{n,k-2}$ follows by taking $q=1$ in part {\bf A.}{\sl (1)} of the theorem statement. 
See that $\mathscr{H}^{\alpha}_{n,k} = n\mathscr{H}^{\alpha}_{n,k-1} - \mathscr{H}^{\alpha}_{n,k-2}$ by taking $q=1$ in part {\bf A.}{\sl (4)}. 
For $n>2$, the formula expressing each of $\mathscr{A}^{\alpha}_{n,k}$, $\mathscr{G}^{\alpha}_{n,k}$, and $\mathscr{H}^{\alpha}_{n,k}$ in terms of the roots of the polynomial $x^{2}-nx+1$ can be obtained by solving the recurrence established in \CountingTheorem\ using, say, standard generating function techniques. 
When $n=2$ and $\alpha=1$, $g(x) := \sum_{k \geq 0}\mathscr{X}^{\alpha}_{n,k}x^{k} = 1/(x-1)^{2} = \frac{d}{dx}(1/(1-x))$, from which the claimed sequence values can be immediately obtained.\hfill\QED

\vspace*{0.1in}
\noindent 
{\large \bf \S \DiscussionSection. Further discussion.} 
The objects and results of this discourse are accessible to undergraduate and beginning graduate students.  
Questions naturally related to this discourse might be useful for student projects. 
Here are some possible examples. 
\begin{itemize}
\item Many interpretations of the Fibonacci and Lucas array numbers can be found in their respective OEIS entries.  
It could be interesting to see how some of these interpretations might naturally extend to the Gibonacci array numbers and/or interact with the order-theoretic aspects of \S \LatticeSection\ above.
\item Outside of the general theory of root geometry as developed in \cite{GMTW1}, we can demonstrate the following results about the outputs of signed-alternating Gibonacci polynomials at the specific input value $x=4$. 
To wit, let $\myq := \alpha/\beta$, and let $m$ be any nonnegative integer. 
We can show that $\polyGqone{2m}(4) = -(2m-1)\myq + 4m$ and $\polyGqone{2m+1}(4) = -(m)\myq + 2m+1$. 
It might be interesting to see if there are other specific inputs whose outputs can be as explicitly understood. 
\item One can view the Gibonaccian Networked-numbers Game on two nodes as requiring one of two initial steps before proceeding with the usual firing  moves of the Networked-numbers Game.  Can this idea be extended to a Gibonaccian game on more than two nodes?
\item All roots of the sign-alternating Fibonacci polynomial $\polyGoneone{2^{n}-1}(x)$ are expressible as certain combinations of the symbols `$\sqrt{\ \ }$', `$2$', and `$\pm$'.  For example, the roots of \[\polyGoneone{15}(x) = x^{7}-14x^{6}+78x^{5}-220x^{4}+330x^{3}-252x^{2}+84x-8\] are, in increasing order, \[\left\{2-\sqrt{2+\sqrt{2}}\, \mathbf{,}\,  2-\sqrt{2}\, \mathbf{,}\,   2-\sqrt{2-\sqrt{2}}\, \mathbf{,}\,   2\, \mathbf{,}\,   2+\sqrt{2-\sqrt{2}}\, \mathbf{,}\,   2+\sqrt{2}\, \mathbf{,}\,   2+\sqrt{2+\sqrt{2}}\right\}.\]  Are there any other sub-families of sign-alternating Gibonacci polynomials where something similar can be said? 
\item The reader who is interested in developing some proficiency with the objects of this paper might consider solving the inclusion-exclusion exercise of Section 5.  Another potentially helpful exercise might be to re-derive the expressions of \ExplicitCorollaryTwo\ for the roots of the sign-alternating Fibonacci and sign-alternating Lucas polynomials from the Binet-type formulas of \ExplicitTheorem. 
\end{itemize}

%
\vspace*{-0.2in}
\renewcommand{\refname}{\Large \bf References}
\renewcommand{\baselinestretch}{1.1}


\small\normalsize
\begin{thebibliography}{999999999} 

\bibitem[BQ]{BQ} A.\ T.\ Benjamin and J.\ J.\ Quinn, {\em Proofs That Count: The Art of Combinatorial Proof}, MAA, 2003.  

\bibitem[Ben]{Benoumhani} M.\ Benoumhani, ``A sequence of binomial coefficients related to Lucas and Fibonacci numbers,'' {\em J. Integer Seq.} {\bf 6} (2003), 10 pp.  

\bibitem[BB]{BB} A.\ Bj\"{o}rner and  F.\ Brenti, {\em Combinatorics of Coxeter Groups}, Springer, New York, 2005. 

\bibitem[D]{DonDiamond} R.\ G.\ Donnelly, ``Diamond-colored modular and distributive lattices'', an {\tt arXiv} manuscript, {\tt arXiv:1812.04434} (2018), 20pp. 

\bibitem[DDMN]{DDMN} R.\ G.\ Donnelly, M.\ W.\ Dunkum, S.\ V.\ Malone, and A.\ Nance, ``Symmetric Fibonaccian lattices and representations of the special linear Lie algebras'', an {\tt arXiv} preprint (2020), 18 pp. 

\bibitem[E]{E} K.\ Eriksson, ``The numbers game and Coxeter groups'', {\em Discrete Math.} {\bf 139} (1995), 155--166.

\bibitem[HK]{HK} T.\ Horzum and E.\ G.\ Kocer, ``On some properties of Horadam polynomials'', {\em International Mathematical Forum} {\bf 4} (2009), 1243--1252

\bibitem[GMTW1]{GMTW1} J.\ L.\ Gross, T.\ Mansour, T.\ W.\ Tucker, and D.\ G.\ L.\ Wang, ``Root geometry of polynomial sequences I: Type $(0,1)$'', {\em Journal of Mathematical Analysis and Applications} {\bf 433} (2016), 1261--1289.

\bibitem[GMTW2]{GMTW2} J.\ L.\ Gross, T.\ Mansour, T.\ W.\ Tucker, and D.\ G.\ L.\ Wang, ``Root geometry of polynomial sequences I: Type $(1,0)$'', {\em Journal of Mathematical Analysis and Applications} {\bf 441} (2016), 499--528.

\bibitem[HB]{HB} V.\ E.\ Hoggatt Jr.\ and M.\ Bicknell, ``Roots of Fibonacci polynomials'', {\em Fibonacci Quarterly} {\bf 11} (1973), 271--274.

\bibitem[Kosh]{Koshy} T.\ Koshy, {\em Fibonacci and Lucas Numbers with Applications}, Vol.\ 1, 2$^{\mbox{\tiny nd}}$ ed., John Wiley \& Sons, Hoboken, New Jersey, 2018.  

\bibitem[KK]{KK} A.\ Khrabrov and A.\ Khokas, ``Points on a line, shoelaces and dominoes'', an {\tt arXiv} manuscript, {\tt arXiv:1505.06309} (2015), 14 pp.

\bibitem[MS]{Munarini--Salvi} E.\ Munarini and N.\ Z.\ Salvi, ``On the rank polynomial of the lattice of order ideals of fences and crowns,''  {\em Discrete Math.} {\bf 259} (2002), 163--177. 

\bibitem[OEIS]{OEIS} OEIS Foundation Inc., {\em The On-Line Encyclopedia of Integer Sequences founded in 1964 by N.\ J.\ A.\ Sloane}, (2020), {\tt http://oeis.org}.
\end{thebibliography}
\end{document}